\documentclass[journal ]{new-aiaa}
\usepackage[utf8]{inputenc}
\usepackage{textcomp}

\usepackage{bm}
\usepackage{subfig}
\usepackage{amssymb}

\usepackage[svgnames]{xcolor,colortbl}

\usepackage[bookmarks=true,pdftex=false,
bookmarksopen=true,pdfborder={0,0,0},linktocpage=true,
colorlinks=true,linkcolor=RoyalBlue,citecolor=RoyalBlue,
urlcolor=RoyalBlue]{hyperref}

\usepackage{graphicx}
\usepackage{amsmath}
\usepackage[version=4]{mhchem}
\usepackage{siunitx}
\usepackage{longtable,tabularx}
\usepackage{acronym}
\acrodef{iod}[IOD]{initial orbit determination}
\acrodef{leo}[LEO]{low Earth orbit}
\acrodef{ads}[ADS]{automatic domain splitting}
\acrodef{da}[DA]{Differential Algebra}
\acrodef{ls}[LS]{least squares}
\acrodef{dim}[DIM]{Doppler integration method}
\acrodef{los}[LOS]{line-of-sight}
\acrodef{rso}[RSO]{resident space objects}
\acrodef{mc}[MC]{Monte Carlo}
\acrodef{taos}[tAOS]{topocentric Attributable Optimized Coordinate System}
\acrodef{tle}[TLE]{Two-Line Elements}

\usepackage{multicol}
\usepackage{multirow}

\newcommand\blfootnote[1]{%
  \begingroup
  \renewcommand\thefootnote{}\footnote{#1}%
  \addtocounter{footnote}{-1}%
  \endgroup
}

\title{Robust initial orbit determination for short-arc Doppler radar observations\blfootnote{\textit{This work has been submitted to the IEEE for possible publication. Copyright may be transferred without notice, after which this version may no longer be accessible.}}\blfootnote{}}

\author{Matteo Losacco \footnote{Post-doc researcher, DCAS,10 Av. Edouard Belin, Toulouse, France}}
\affil{ISAE-SUPAERO, Toulouse, 31400, France}
\author{Roberto Armellin\footnote{Professor, Te P\=unaha \=Atea - Space Institute, 20 Symonds Street, Auckland, New Zealand}}
\affil{University of Auckland, Auckland, 1010, New Zealand}
\author{Carlos Yanez \footnote{Space Surveillance Engineer, 18 Av. Edouard Belin, Toulouse, France}}
\affil{CNES, Toulouse, 31400, France}
\author{Stéphanie Lizy-Destrez \footnote{Associate Professor in Space Systems, DCAS,10 Av. Edouard Belin, Toulouse, France}}
\affil{ISAE-SUPAERO, Toulouse, 31400, France}
\author{Laura Pirovano\footnote{Post-doc researcher, Te P\=unaha \=Atea - Space Institute, 20 Symonds Street, Auckland, New Zealand}}
\affil{University of Auckland, Auckland, 1010, New Zealand}
\author{Francesco Sanfedino \footnote{Associate Professor in Control of Space Systems, DCAS,10 Av. Edouard Belin, Toulouse, France}}
\affil{ISAE-SUPAERO, Toulouse, 31400, France}

\begin{document}

\maketitle

\begin{abstract}
A new Doppler radar initial orbit determination algorithm with embedded uncertainty quantification capabilities is presented. The method is based on a combination of Gauss' and Lambert's solvers. The whole process is carried out in the Differential Algebra framework, which provides the Taylor expansion of the state estimate with respect to the measurements' uncertainties. This feature makes the approach particularly suited for handling data association problems. A comparison with the Doppler integration method is performed using both simulated and real data. The proposed approach is shown to be more accurate and robust, and particularly suited for short-arc observations.  
\end{abstract}

\section{Introduction}
The accurate estimation of the state of objects orbiting around the Earth is of great importance whenever operations such as observation scheduling, data association and collision risk assessment are required. The quality of these estimates is mostly affected by the accuracy of the sensors used for the observation and the frequency of these observations. As a result, several networks of ground-based optical~\citep{Schildknecht2007}, radar~\citep{Muntoni2021}, and laser~\citep{Cordelli2020} sensors operate all over the world to continuously monitor the near-Earth environment~\citep{Faucher2019}. Nevertheless, the available catalogues of \ac{rso} cover just a small portion of the existing population. The increasingly growing launch activity and in-orbit object generation events, such as collisions, fragmentations and explosions, represent a continuous source of new objects~\citep{ESASDO2021}. Though the size and properties of most of these objects prevent existing operational sensors from detecting or observing them with enough accuracy, observation campaigns often offer the possibility of identifying and potentially characterising uncatalogued objects. 

The process of state estimation of an uncatalogued object starting from a set of  measurements is known in literature as \ac{iod}. Such measurements are typically referred to as tracks, tracklets, or short-arcs. Conversely, passages that are too limited in time to perform \ac{iod} are generally labelled as too-short-arcs, and shall be processed by relying on different approaches, such as the admissible region theory~\citep{Milani2004,DeMars2012,Pirovano2021}. This paper focuses on the first scenario, thus investigates the use of \ac{iod} methods for tracklets and short-arcs.

Classical \ac{iod} algorithms can be divided according to the type of measurements and sensors. Two main families of methods exist, i.e. angle-only methods for optical sensors and angle-range methods for range radars. The first solutions for the angle-only \ac{iod} problem were developed in the eighteenth century, aiming at determining the orbit of celestial bodies such as planets or asteroids. Over more than 200 years, several techniques have been proposed. Among these, Laplace's method~\citep{Laplace1780}, Gauss' method~\citep{Gauss1857}, the double R iteration method~\citep{Escobal1965}, Baker-Jacobi's method~\citep{BakerJr1977}, Gooding's method~\citep{Gooding1996} and Karimi and Mortari's methods~\citep{Karimi2010} represent the most common ones. 
The second family of \ac{iod} methods processes range radars observables, i.e. angular and range measurements. Examples of these so-called angle-range IOD methods are Lambert's method~\citep{Lambert1761,Izzo2015}, Gibbs' method~\citep{Gibbs1889}, and Herrick-Gibbs' method~\citep{Herrick1971}, along with recent approaches tailored for short-arc observations~\citep{Zhang2020}. 

While a vast literature for range radars \ac{iod} exists, few algorithms are available for Doppler radars. These sensors measure the angular position of the object and the time derivative of its slant range $d$~\citep{Barlow1949}, which is defined as the sum of the distances from the radar receiver and transmitter, both in case they coincide (monostatic radar) or not (bistatic configuration). Algorithms such as the \ac{dim}~\citep{Yanez2017} and the hodograph method~\citep{Christian2021} have been recently proposed, proving to provide accurate \ac{iod} estimates when processing Doppler radar observables. Nevertheless, none of these methods was explicitly tested against the challenging case of short-arc observations, which can be common when targeting the \ac{leo} regime.

Another main limitation of classical \ac{iod} methods is the lack of information on the estimate uncertainty, unless performing \ac{ls} refinement or making model simplifications \citep{Zhang2020}. The estimation of the IOD solution uncertainty is however essential to establish robust data association schemes enabling catalogue initialization. To overcome these limitations, new methods based on the use of \ac{da}~\citep{Valli2013} have been recently presented. Armellin and Di Lizia~\citep{Armellin2018} reformulated Lambert's problem in \ac{da} sense, providing a mathematical description of the solution uncertainty of the angle-range IOD problem as a function of measurements noise.  This result, combined with the use of \ac{ads}~\citep{Wittig2015}, was then exploited in \cite{Pirovano2020b} to solve the angle-only \ac{iod} problem, thus allowing for an analytical description of the so-called \textit{orbit set}. 

Starting from these considerations, this paper introduces a new robust \ac{iod} algorithm for Doppler radars with built-in uncertainty quantification capabilities. Given a set of angular and range rate measurements, the method  provides an accurate solution to the \ac{iod} problem even for short-arc observations. The solution uncertainty quantification is achieved by combining \ac{da} and \ac{ads}, which provide the Taylor expansion of the state with respect to the observations. The algorithm, here referred to as DAIOD algorithm, represents an extension of its optical~\citep{Pirovano2020b} and range radar~\citep{Armellin2018} counterparts, and completes the DAIOD estimation framework for the three possible available scenarios: angle-only, angle-range, and angle-Doppler. 

The paper is structured as follows. Section~\ref{sec:DA} describes the required mathematical tools, i.e. \ac{da} and \ac{ads}. A formulation of the method is given in Section~\ref{sec:DAIOD}. Finally, Section~\ref{sec:performance} illustrates the performance of the algorithm, and compares it with the \ac{dim} method.

\section{Differential algebra and Automatic Domain Splitting}
\label{sec:DA}
Differential algebra provides the tools to compute the derivatives of functions in a computer environment~\citep{Ritt1932,Kolchin1973}. The basic idea of \ac{da} is to bring the treatment of functions and the operations on them to the computer environment in a similar way as the treatment of real numbers. More specifically, by substituting the algebra of real numbers with a new algebra of Taylor polynomials, any sufficiently regular function $\bm{f}$ of $v$ variables can be expanded into its Taylor expansion up to an arbitrary order $k$. Starting from this expression, the Taylor coefficients of order $k$ representing sum and product of functions, or scalar products with real numbers, can be directly computed from the coefficients of their summands and factors. As a consequence, the set of equivalence classes of functions can be endowed with well-defined operations, leading to the so-called truncated power series algebra. In addition to basic algebraic operations, differentiation and integration can be easily introduced in the algebra, thus finalizing the definition of the differential algebra structure of DA~\citep{Berz1999}. 

A detailed formulation of DA for astrodynamics applications can be found in~\cite{DiLizia2008}. As an example, consider a generic multivariate function $\bm{y}=\bm{f}(\bm{x})$, $\bm{x}\in\mathbb{R}^v$. Starting from a nominal $\bm{x}$ and the associated uncertainty $\Delta$, here assumed equal for all components, the \ac{da} representation of $\bm{x}$ can be expressed as
\begin{equation}
\label{eq:x_DA}
[\bm{x}]=\bm{x}+\Delta\delta\bm{x}
\end{equation}
with $\delta\bm{x}\in\mathbb{R}^v$ representing the deviation from $\bm{x}$. If $\bm{f}$ is evaluated in the \ac{da} framework, one obtains
\begin{equation}
\label{eq:y_DA}
[\bm{y}] = \bm{f}([\bm{x}])=\mathcal{T}_{\bm{y}}(\delta\bm{x})
\end{equation}
The term $\mathcal{T}_{\bm{y}}(\delta\bm{x})$ of Eq.~\eqref{eq:y_DA} indicates the Taylor expansion of $\bm{y}$ with respect to $\delta\bm{x}$. As a result, by considering a generic initial deviation $\delta\bm{x}^*$, the corresponding $\bm{y}^*$ solution can be directly obtained by mapping $\delta\bm{x}^*$ with $\mathcal{T}_{\bm{y}}(\delta\bm{x})$, that is
\begin{equation}
\label{eq:y_star}
\bm{y}^*=\mathcal{T}_{\bm{y}}(\delta\bm{x}^*)
\end{equation}

The combined use of \ac{da} and polynomial bounding techniques~\citep{Crane1975} offers a powerful tool to estimate the bound associated with $[\bm{y}]$. Consider a generic component $j$ of $\bm{y}$, and define an interval for $\delta\bm{x}$, e.g. $\delta\bm{x}\in[-1,1]^{v}$. The use of polynomial bounders allows us to write
\begin{equation}
y_j \in \left[-b_{y_j}, b_{y_j}\right]
\end{equation}
This functionality will be exploited in the paper to estimate the uncertainty of the computed \ac{iod} solution.

\subsection{The \ac{ads} algorithm}
\label{subsec:ADS}
The accuracy of the \ac{da} expansion tends to decrease for an increase in the size of the uncertainty set $\delta\bm{x}$ and for higher nonlinearity levels of $\bm{f}$. A possible approach would consist in properly increasing the selected expansion order $k$. This, however, does not always improve the expansion accuracy, while typically implies an increase in the required computational burden, which may make the \ac{da} formulation unfeasible in case of highly nonlinear functions.  On the other hand, the order is not the only available parameter to play with for increasing the accuracy. Suppose to maintain the order unaltered and to split the whole uncertainty set into smaller subsets. By dividing the initial uncertainty domain and computing the Taylor expansion around the center points of the new sets, the same overall coverage is granted, but with a reduced set size (and thus a larger accuracy) per expansion. Starting from these considerations, \ac{ads} employs an automatic algorithm to determine whether the current polynomial representation is accurate enough or not~\citep{Wittig2015}. If the accuracy requirements are not met, the domain of the original expansion is divided along one of the expansion variables into two domains of half their original size. By re-expanding the polynomials around the new centre points, two separate expansions are obtained, each covering a half of the original set. At this point, the procedure is repeated on the generated sets, and it terminates as soon as the desired accuracy is obtained or a threshold on the number of splits is reached. At the end of the process, a so-called \textit{manifold} of Taylor expansions is obtained. If we recall the previous example, the solution provided by the \ac{ads} can be written as
\begin{equation}
[\bm{y}]=\bigcup\limits_{i=1}^{N_s}\mathcal{T}_{\bm{y}}^i(\delta\bm{x})
\end{equation}
where $N_s$ is the number of generated sets. Two parameters govern the accuracy of the \ac{ads} result, namely the tolerance for the splitting decision and the maximum number of splits per direction. A detailed description of the method can be found in~\cite{Wittig2015}.

\section{The DAIOD method}
\label{sec:DAIOD}
Consider a set of angular and range rate measurements provided by a Doppler radar
\begin{equation}
\label{eq:meas_topo_set}
\begin{aligned}
\left(t_i;\left(A_i;\sigma_{A_i}\right),\left(E_i;\sigma_{E_i}\right),\left(\dot{d}_i;\sigma_{\dot{d}_i}\right)\right)
\end{aligned}
\end{equation}
with $i\in[1,N]$, while $A_i$, $E_i$ and $\dot{d}_i$ are the azimuth, elevation, and range rate measurements obtained at epoch $t_i$, respectively, whereas $\sigma_{A_i}$, $\sigma_{E_i}$ and $\sigma_{\dot{d}_i}$ are the associated standard deviations of the sensor noise. The DAIOD algorithm estimates the state of the orbiting object at the epoch of the first available set of measurements by processing angular and range rate data. The method works in three phases: measurements processing, \ac{iod} solution computation, and \ac{iod} solution expansion.

\subsection{Phase 1: Measurements processing}
\label{subsec:DAIOD_phase1}
Sensor measurements can be typically treated as independent Gaussian random variables. As a result, each measurement is characterised by the mean value and the standard deviation. Therefore, Eq.~\eqref{eq:meas_topo_set} can be re-expressed as
\begin{equation}
\label{eq:meas_topo_rand}
\begin{aligned}
\left(t_i;Y_{A_i},Y_{E_i},Y_{\dot{d}_i}\right)
\end{aligned}
\end{equation}
where $Y_{\beta_i}\sim\mathcal{N}\left(\beta_i,\sigma_{\beta_i}\right)$ indicates a normal random distribution for the generic measurement $\beta_i$.

The DAIOD algorithm requires as input two angles and a value of range rate per time instant. In the adopted formulation, the two angles are the topocentric right ascension $\alpha$ and declination $\delta$. Two different options are available. One may choose to directly use the available raw data to perform the \ac{iod} estimation process. In this case, the angular measurements of Eq.~\eqref{eq:meas_topo_rand} are converted into $\alpha$ and $\delta$ by using an uncertainty mapping method. For the case under study, simple \ac{mc} is used. More specifically, given a time epoch $t_i$ and the associated $Y_{A_i}$ and $Y_{E_i}$, a set of $N_{MC}$ samples $\{A_i^j;E_i^j\}$ is generated. Each couple is converted into the corresponding $\{\alpha_i^j;\delta_i^j\}$. At this point, the estimated ($\hat{\quad}$ superscript) mean values $\hat{\alpha}_i$ and $\hat{\delta}_i$ can be easily computed. In addition, given a desired confidence level $\iota_L\in[0,1]$, the corresponding confidence interval amplitude can be derived as
\begin{equation}
\Delta \textrm{CI}_{\hat{\beta}_i}(\iota_L) = \left\{\textrm{d}\left(\beta_i^{j^*},\hat{\beta}_i\right) : p\left(\textrm{d}\left(\beta_i^{j},\hat{\beta}_i\right)<\textrm{d}\left(\beta_i^{j^*},\hat{\beta}_i\right)\right)\geq\iota_L\right\}
\end{equation}
with $i\in[1,N]$, $j\in[1,N_{MC}]$, whereas $\textrm{d}\left(a,b\right)$ refers to the difference between angles $a$ and $b$ and $p$ is the probability. As a result, for each time instant we can write
\begin{equation}
\label{eq:alpha_delta}
\begin{gathered}
\textrm{CI}_{\hat{\alpha}_i} = \left[\hat{\alpha}_i\pm\Delta\textrm{CI}_{\hat{\alpha}_i}(\iota_L)\right]\\
\textrm{CI}_{\hat{\delta}_i} = \left[\hat{\delta}_i\pm\Delta\textrm{CI}_{\hat{\delta}_i}(\iota_L)\right]\\
\end{gathered}
\end{equation} 
The confidence interval for range rate measurements is instead directly derived from the sensor accuracy, i.e.
\begin{equation}
\textrm{CI}_{\hat{\dot{d}}_i} = \left[\hat{\dot{d}}_i\pm k(\iota_L)\sigma_{\hat{\dot{d}}_i}\right]
\end{equation}
where $\hat{\dot{d}}_i = \dot{d}_i$, $\sigma_{\hat{\dot{d}}_i} = \sigma_{\dot{d}_i}$ while $k(\iota_L)$ is a scaling of the noise standard deviation.

The second available option is to perform regression on the measurements to reduce the effect of the noise. In this case, the choice of the reference frame in which the regression is performed is crucial to limit the dependency of the regression order on passage duration and shape. A detailed analysis for range and Doppler radars is illustrated in~\cite{Reihs2021}. The reference frame adopted in this work is the so-called \ac{taos}. The centre of this frame coincides with the origin of the topocentric frame of the radar receiver, and the angular position of the space object is computed with respect to the plane defined by the first and last \ac{los} of the passage. As a result, a generic \ac{los} can be expressed in terms of two angles, here defined $\lambda$ and $\gamma$. The first step of the regression process consists in converting all the $\left\{A_i;E_i\right\}$ couples into their $\left\{\lambda_i;\gamma_i\right\}$ counterpart. The procedure is the same as the one for raw data (the $\hat{\,}$ superscript is here dropped for convenience). Range rate measurements are instead unaltered, since the topocentric and the \ac{taos} frames share the same origin and they are obtained with a time-independent rotation. Once defined the new set of measurements, the regression problem can be setup. By defining with $\beta$ the generic observable, $\beta=\{\lambda,\gamma,\dot{d}\}$, the following system can be written
\begin{equation}
\begin{bmatrix}
\beta_1\\
\beta_2\\
\vdots\\
\beta_N
\end{bmatrix}=\bm{A}_{\beta}\bm{z}_{\beta}
\end{equation}
The explicit expressions for the coefficient matrix $\bm{A}_{\beta}$ and the vector of unknowns $\bm{z}_{\beta}$ depend on the selected regression order. For a generic order $n$, we can write
\begin{equation}
\label{eq:A}
\bm{A}_{\beta}=\begin{bmatrix}
1& t_1-t_0& \cdots& (t_1-t_0)^n\\
1& t_2-t_0& \cdots& (t_2-t_0)^n\\
\vdots& \vdots& \cdots& \vdots\\
1& t_N-t_0& \cdots& (t_N-t_0)^n
\end{bmatrix}
\end{equation}
while $\bm{z}_{\beta}=\left[\beta_0\,\dot{\beta}_0\cdots\overset{n}{\dot{\beta}}_0\right]^{\textrm{T}}$.
In Eq.~\eqref{eq:A}, $t_0$ indicates the regression epoch, which is selected equal to the mid-observation time. Once the regression order is selected, the problem becomes
\begin{equation}
\bm{y}_{\beta}=\bm{A}_{\beta}\bm{z}_{\beta}
\end{equation}
The solution of the problem follows the classical \ac{ls} theory. The normal matrix $\bm{N}_{\beta}$ can be written as
\begin{equation}
\label{eq:Regr_N}
\bm{N}_{\beta}=\bm{A}_{\beta}^{\textrm{T}}\bm{Q}_{\beta}^{-1}\bm{A}_{\beta}=\bm{A}_{\beta}^{\textrm{T}}\bm{A}_{\beta}
\end{equation}
where the cofactor matrix $\bm{Q}_{\beta} = \bm{I}$ assuming zero mean, single variance and independent errors. As a result, an estimate of the vector of unknown parameters can be obtained
\begin{equation}
\hat{\bm{z}}_{\beta}=\bm{N}_{\beta}^{-1}\bm{A}_{\beta}^{\textrm{T}}\bm{y}_{\beta}
\end{equation}
The availability of $\hat{\bm{z}}_{\beta}$ allows us to compute an estimate of the observations and the associated residuals
\begin{equation}
\begin{aligned}
\hat{\bm{y}}_{\beta}&=\bm{A}_{\beta}\hat{\bm{z}}_{\beta}\\
\hat{\bm{c}}_{\beta} &= \bm{y}_{\beta}-\hat{\bm{y}}_{\beta}
\end{aligned}
\end{equation}
As a result, an estimate for the standard deviation of the observations can be obtained
\begin{equation}
\hat{\sigma}^2_{\beta}=\dfrac{\hat{\bm{c}}_{\beta}^{\textrm{T}}\bm{Q}_{\beta}^{-1}\hat{\bm{c}}_{\beta}}{N-p}=\dfrac{\hat{\bm{c}}_{\beta}^{\textrm{T}}\hat{\bm{c}}_{\beta}}{N-p}
\end{equation}
where $p=n+1$. This estimate is finally used to compute the covariance matrix of the unknown parameters $\bm{C}_{\hat{\bm{z}}_{\beta}\hat{\bm{z}}_{\beta}}$ and the covariance of the estimated observations $\bm{C}_{\hat{\bm{y}}_{\beta}\hat{\bm{y}}_{\beta}}$
\begin{equation}
\label{eq:Regr_cov}
\begin{aligned}
\bm{C}_{\hat{\bm{z}}_{\beta}\hat{\bm{z}}_{\beta}}&=\hat{\sigma}^2_{\beta}\bm{N}_{\beta}^{-1}\\
\bm{C}_{\hat{\bm{y}}_{\beta}\hat{\bm{y}}_{\beta}} &=\bm{A}_{\beta}\bm{C}_{\hat{\bm{z}}_{\beta}\hat{\bm{z}}_{\beta}}\bm{A}_{\beta}^{\textrm{T}}
\end{aligned}
\end{equation}
Let us now assume that the measurements are distributed as independent normal variables. This is valid for $\dot{d}$, but it is just an approximation for the tAOS angles $\lambda$ and $\gamma$, as the uncertainty mapping from $A$ and $E$ does not guarantee the Gaussianity conservation. If this holds, it can be shown that
\begin{equation}
\dfrac{\hat{\beta}_i-\beta_i}{\sqrt{C_{\hat{\bm{y}}_{\beta}\hat{\bm{y}}_{\beta},ii}}} \sim t_{N-p}
\end{equation}
where $t_{N-p}$ is the Student's t-distribution with $N-p$ degrees of freedom, whereas $C_{\hat{\bm{y}}_{\beta}\hat{\bm{y}}_{\beta},ii}$ is the $i$-th diagonal element of matrix $\bm{C}_{\hat{\bm{y}}_{\beta}\hat{\bm{y}}_{\beta}}$, $i\in[1,N]$. As a result, given the desired confidence level $\iota_L$  and the associated $t_{N-p}$ quantile $q_{\iota_L}$, the confidence intervals can be written as
\begin{equation}
\label{eq:CI_tAOS}
\textrm{CI}_{\hat{\beta}_i} = \left[\hat{y}_{\beta,i}\pm q_{\iota_L}\sqrt{C_{\hat{\bm{y}}_{\beta}\hat{\bm{y}}_{\beta},ii}}\right]
\end{equation}
Equation~\eqref{eq:CI_tAOS} provides a statistical representation for $(\hat{\lambda}_i;\hat{\gamma}_i;\hat{\dot{d}}_i)$ at all time epochs. At this point, the angular measurements $\hat{\lambda}_i$ and $\hat{\gamma}_i$ are converted into $\hat{\alpha}_i$ and $\hat{\delta}_i$ through a procedure similar to the one described for raw data, thus obtaining the formulation of Eq.~\eqref{eq:alpha_delta}.

At the end of phase 1, regardless of the selected measurement processing, the following set is available
\begin{equation}
\label{eq:meas_polar_set}
\begin{aligned}
\left(t_i;\left(\hat{\alpha}_i;\Delta\textrm{CI}_{\hat{\alpha}_i}\right),\left(\hat{\delta}_i;\Delta\textrm{CI}_{\hat{\delta}_i}\right),\left(\hat{\dot{d}}_i;\Delta\textrm{CI}_{\hat{\dot{d}}_i}\right)\right)
\end{aligned}
\end{equation}

\subsection{Phase 2: \ac{iod} solution computation}
\label{subsec:DAIOD_phase2}
The second phase of the DAIOD method combines a Gauss and an iterative Lambert's solver to estimate the object range $\rho$ from the radar receiver at the epochs of the first and last measurements. To limit the impact of Gauss' solver failures, which are likely to occur with short-arc observations~\citep{Pirovano2020b}, we adopt here a strategy that leverages the intrinsically large uncertainty of radar angular measurements to grant high success rates regardless of the arc length. The algorithm exploits DA and map inversion techniques~\citep{Berz1999} to solve the resulting system of implicit equations, thus avoiding additional computations typically required by standard nonlinear solvers, e.g. Newton-Raphson method.

Let us start from Eq.~\eqref{eq:meas_polar_set} and restrict the analysis to the angular measurements $\hat{\alpha}$ and $\hat{\delta}$ at the first (1), middle ($m$), and last ($N$) observation epochs. If one takes their nominal value plus all possible combinations built by considering the two extremes of their confidence intervals, 65 different sets are obtained, i.e.
\begin{equation}
\label{eq:sets}
\begin{aligned}
\left(\hat{\alpha}_1^j,\hat{\delta}_1^j,\hat{\alpha}_m^j,\hat{\delta}_m^j,\hat{\alpha}_N^j,\hat{\delta}_N^j\right)
\end{aligned}\qquad j\in[1,65]
\end{equation}
For each set of Eq.~\eqref{eq:sets}, this procedure is followed:
\begin{enumerate}
\item[\label=\color{RoyalBlue} 1.] Solve the Gauss problem, $\mathcal{G}$, thus obtaining an estimate for the range at $t_1$, $t_m$, and $t_N$
\begin{equation}
\begin{gathered}
\mathcal{G}\left(t_1,\hat{\alpha}_1^j,\hat{\delta}_1^j,t_m,\hat{\alpha}_m^j,\hat{\delta}_m^j,t_N,\hat{\alpha}_N^j,\hat{\delta}_N^j\right)\rightarrow\\\left\{\hat{\rho}_{1,\mathcal{G}}^{j},\hat{\rho}_{m,\mathcal{G}}^{j},\hat{\rho}_{N,\mathcal{G}}^{j}\right\} 
\end{gathered}
\end{equation}
\item[\label=\color{RoyalBlue} 2.] Set $\hat{\rho}_1^{k,j}=\hat{\rho}_{1,\mathcal{G}}^{j}$ and $\hat{\rho}_N^{k,j}=\hat{\rho}_{N,\mathcal{G}}^{j}$, with $k=1$, and initialize the ranges as DA variables
\begin{equation}
\begin{aligned}
    \left[\hat{\rho}_1^{k,j}\right] &=\hat{\rho}_1^{k,j}+\delta\rho_1\\
    \left[\hat{\rho}_N^{k,j}\right] &=\hat{\rho}_N^{k,j}+\delta\rho_N
\end{aligned}
\end{equation}

\item[\label=\color{RoyalBlue} 3.] Solve the DA Lambert's problem~\citep{Armellin2018}
\begin{equation}
\begin{gathered}
\mathcal{L}\left(t_1,\hat{\alpha}_1,\hat{\delta}_1,\left[\hat{\rho}_1^{k,j}\right], t_N, \hat{\alpha}_N,\hat{\delta}_N,\left[\hat{\rho}_N^{k,j}\right]\right)\rightarrow\\
\begin{Bmatrix}
\left[\hat{\bm{v}}_{1,\mathcal{L}}^{k,j}\right]=\mathcal{T}_{\hat{\bm{v}}_{1,\mathcal{L}}^{k,j}}(\delta\bm{\rho})
\vspace{0.2cm}\\
\left[\hat{\bm{v}}_{N,\mathcal{L}}^{k,j}\right]=\mathcal{T}_{\hat{\bm{v}}_{N,\mathcal{L}}^{k,j}}(\delta\bm{\rho})
\end{Bmatrix}
\end{gathered}
\end{equation}
where $\delta\bm{\rho} = \left\{\delta\rho_1,\delta\rho_N\right\}$, while $\bm{v}$ is the velocity.
\item[\label=\color{RoyalBlue} 4.] Compute the expansion of the range rate at $t_1$ and $t_N$
\begin{equation}
\begin{aligned}
\left[\hat{\dot{d}}_{1,\mathcal{L}}^{k,j}\right] &= \mathcal{T}_{\hat{\dot{d}}_{1,\mathcal{L}}^{k,j}}(\delta\bm{\rho})\\
\left[\hat{\dot{d}}_{N,\mathcal{L}}^{k,j}\right] &= \mathcal{T}_{\hat{\dot{d}}_{N,\mathcal{L}}^{k,j}}(\delta\bm{\rho})
\end{aligned}
\end{equation}
\item[\label=\color{RoyalBlue} 5.] Compute the deviations with respect to $\hat{\dot{d}}_1$ and $\hat{\dot{d}}_N$
\begin{equation}
\left[\Delta\hat{\dot{\bm{d}}}^{k,j}\right]=
\begin{Bmatrix}
\left[\hat{\dot{d}}_{1,\mathcal{L}}^{k,j}\right]-\hat{\dot{d}}_1
\vspace{0.2cm}\\
\left[\hat{\dot{d}}_{N,\mathcal{L}}^{k,j}\right]-\hat{\dot{d}}_N
\end{Bmatrix}
\end{equation}
\item[\label=\color{RoyalBlue} 6.] Consider the origin preserving map
\begin{equation}
\delta\dot{\bm{d}} = \left[\Delta\hat{\dot{\bm{d}}}^{k,j}\right]-\Delta\hat{\dot{\bm{d}}}^{k,j} = \mathcal{T}_{\delta\dot{\bm{d}}}\left(\delta\bm{\rho}\right)
\end{equation}
\item[\label=\color{RoyalBlue} 7.] Invert the map
\begin{equation}
\label{eq:inv}
\delta\bm{\rho} = \mathcal{T}_{\delta\bm{\rho}}\left(\delta\dot{\bm{d}}\right)
\end{equation}
\item[\label=\color{RoyalBlue} 8.] Evaluate Eq.~\eqref{eq:inv} in $-\Delta\hat{\dot{\bm{d}}}^{k,j}$
\begin{equation}
\Delta\hat{\bm{\rho}}^{k,j} = \mathcal{T}_{\delta\bm{\rho}}\left(-\Delta\hat{\dot{\bm{d}}}^{k,j}\right)=
\begin{Bmatrix}
\Delta\hat{\rho}_1^{k,j}\\
\Delta\hat{\rho}_N^{k,j}
\end{Bmatrix}
\end{equation}
\item[\label=\color{RoyalBlue} 9.] Update the ranges
\begin{equation}
\begin{aligned}
\left[\hat{\rho}_1^{k+1,j}\right] &= [\hat{\rho}_1^{k,j}]+\Delta\hat{\rho}_1^{k,j}\\
[\hat{\rho}_N^{k+1,j}] &= [\hat{\rho}_N^{k,j}]+\Delta\hat{\rho}_N^{k,j}
\end{aligned}
\end{equation}
\item[\label=\color{RoyalBlue} 10.] Check $\lvert\lvert\hat{\bm{\rho}}^{k+1,j}-\hat{\bm{\rho}}^{k,j}\rvert\rvert$: if it is lower than an imposed threshold stop, otherwise set $k=k+1$ and return to step 3.
\end{enumerate}
The described process allows us to obtain multiple estimates for the range values at $t_1$ and $t_N$, one for each sample $j$, and terminates as soon as two solutions $(\hat{\rho}_1^k;\hat{\rho}_N^k)$ and $(\hat{\rho}_1^h;\hat{\rho}_N^h)$ whose difference is below a given threshold are found, or all the 65 sets of angles are investigated. Then, the most likely solution is identified by considering the residual with respect to the whole set of available measurements. That is, knowing
\begin{equation}
 \left(\hat{\alpha}_1,\hat{\delta}_1,\hat{\rho}_1^j, \hat{\bm{v}}_{1,\mathcal{L}}^j\right)\rightarrow\left(\hat{\bm{r}}_{1}^j,\hat{\bm{v}}_{1,\mathcal{L}}^j\right)  
\end{equation}
where $\bm{r}$ is the object position vector, an estimate for $\alpha$, $\delta$ and $\dot{d}$ at all $t_i$ can be obtained, i.e.
\begin{equation}
\label{eq:meas_polar_lamb}
\left(\hat{\alpha}^j_{i,\mathcal{L}},\hat{\delta}^j_{i,\mathcal{L}},\hat{\dot{d}}^j_{i,\mathcal{L}}\right)\qquad i\in[1,N]
\end{equation}
The residual between the measurements given by Eq.~\eqref{eq:meas_polar_set} and those given by Eq.~\eqref{eq:meas_polar_lamb} can be expressed as
\begin{equation}
R^j = \sum_{i=1}^N  \left(\left(\dfrac{\hat{\alpha}_i-\hat{\alpha}^j_{i,\mathcal{L}}}{\Delta\textrm{CI}_{\hat{\alpha}_i}}\right)^2+\left(\dfrac{\hat{\delta}_i-\hat{\delta}^j_{i,\mathcal{L}}}{\Delta\textrm{CI}_{\hat{\delta}_i}}\right)^2+\left(\dfrac{\hat{\dot{d}}_i-\hat{\dot{d}}^j_{i,\mathcal{L}}}{\Delta\textrm{CI}_{\hat{\dot{d}}_i}}\right)^2\right)
\end{equation}
The selected $(\hat{\rho}_1^{j^*};\hat{\rho}_N^{j^*})$ is the one with the lowest $R^j$.

\subsection{Phase 3: \ac{iod} solution expansion}
\label{subsec:DAIOD_phase3}
At the end of phase 2 an estimate for the \ac{rso} state at the epoch of the first measurement is obtained. The third phase of the DAIOD algorithm exploits the DA Lambert's solver to expand the \ac{iod} solution in Taylor series with respect to the observables. The steps are the following:
\begin{enumerate}
\item[\label=\color{RoyalBlue} 1.] Initialize the angles and range rate measurements at epochs 1 and $N$ as \ac{da} variables
\begin{equation}
\begin{aligned}
\left[\hat{\alpha}_1\right] &= \hat{\alpha}_1+\Delta\textrm{CI}_{\hat{\alpha}_1}\delta\alpha_1\\
\left[\hat{\delta}_1\right] &= \hat{\delta}_1+\Delta\textrm{CI}_{\hat{\delta}_1}\delta\delta_1\\
\left[\hat{\dot{d}}_1\right] &= \hat{\dot{d}}_1+\Delta\textrm{CI}_{\hat{\dot{d}}_1}\delta\dot{d}_1\\
\left[\hat{\alpha}_N\right] &= \hat{\alpha}_N+\Delta\textrm{CI}_{\hat{\alpha}_N}\delta\alpha_N\\
\left[\hat{\delta}_N\right] &= \hat{\delta}_N+\Delta\textrm{CI}_{\hat{\delta}_N}\delta\delta_N\\
\left[\hat{\dot{d}}_N\right] &= \hat{\dot{d}}_N+\Delta\textrm{CI}_{\hat{\dot{d}}_N}\delta\dot{d}_N
\end{aligned}
\end{equation}
\item[\label=\color{RoyalBlue} 2.] Initialize the ranges as \ac{da} variables
\begin{equation}
\begin{gathered}
\left[\hat{\rho}_1\right] = \hat{\rho}_1^{j^*}+\delta\rho_1\qquad
\left[\hat{\rho}_N\right] = \hat{\rho}_N^{j^*}+\delta\rho_N
\end{gathered}
\end{equation}
\item[\label=\color{RoyalBlue} 3.] Solve the associated DA Lambert's problem
\begin{equation}
\begin{gathered}
\mathcal{L}\left(t_1,\left[\hat{\alpha}_1\right],\left[\hat{\delta}_1\right],\left[\hat{\rho}_1\right], t_N,\left[\hat{\alpha}_N\right],\left[\hat{\delta}_N\right],\left[\hat{\rho}_N\right]\right)\rightarrow\\
\begin{Bmatrix}
\left[\hat{\bm{v}}_{1,\mathcal{L}}\right]=\mathcal{T}_{\hat{\bm{v}}_{1,\mathcal{L}}}(\delta\bm{\alpha},\delta\bm{\delta},\delta\bm{\rho})
\vspace{0.2cm}\\
\left[\hat{\bm{v}}_{N,\mathcal{L}}\right]=\mathcal{T}_{\hat{\bm{v}}_{N,\mathcal{L}}}(\delta\bm{\alpha},\delta\bm{\delta},\delta\bm{\rho})
\end{Bmatrix}
\end{gathered}
\end{equation}
where $\delta\bm{\alpha}=\left\{\delta\alpha_1,\delta\alpha_N\right\}$, $\delta\bm{\delta}=\left\{\delta\delta_1,\delta\delta_N\right\}$ and $\delta\bm{\rho}=\left\{\delta\rho_1,\delta\rho_N\right\}$

\item[\label=\color{RoyalBlue} 4.] Compute the Taylor expansion of the range rates
\begin{equation}
\begin{aligned}
\left[\hat{\dot{d}}_{1,\mathcal{L}}\right] &= \mathcal{T}_{\hat{\dot{d}}_{1,\mathcal{L}}}(\delta\bm{\alpha},\delta\bm{\delta},\delta\bm{\rho})\\
\left[\hat{\dot{d}}_{N,\mathcal{L}}\right] &= \mathcal{T}_{\hat{\dot{d}}_{N,\mathcal{L}}}(\delta\bm{\alpha},\delta\bm{\delta},\delta\bm{\rho})
\end{aligned}
\end{equation}
\item[\label=\color{RoyalBlue} 5.] Compute the deviations with respect to $\hat{\dot{d}}_1$ and $\hat{\dot{d}}_N$
\begin{equation}\label{eq:map_partial}
\begin{gathered}
\delta\dot{\bm{d}}=
\begin{Bmatrix}
\left[\hat{\dot{d}}_{1,\mathcal{L}}\right]-\hat{\dot{d}}_1
\vspace{0.2cm}\\
\left[\hat{\dot{d}}_{N,\mathcal{L}}\right]-\hat{\dot{d}}_N
\end{Bmatrix} =
\mathcal{T}_{\delta\dot{\bm{d}}}(\delta\bm{\alpha},\delta\bm{\delta},\delta\bm{\rho})
\end{gathered}
\end{equation}
\item[\label=\color{RoyalBlue} 6.] Append the identities
\begin{equation}
\begin{Bmatrix}
\delta\bm{\alpha}\\
\delta\bm{\delta}\\
\delta\dot{\bm{d}}\\
\end{Bmatrix}
=
\begin{Bmatrix}
\delta\bm{\alpha}\\
\delta\bm{\delta}\\
\mathcal{T}_{\delta\dot{\bm{d}}}(\delta\bm{\alpha},\delta\bm{\delta},\delta\bm{\rho})
\end{Bmatrix}
\end{equation}
\item[\label=\color{RoyalBlue} 7.] Invert the map, and compose it with the available range rate measurements map
\begin{equation}
\label{eq:map}
\delta\bm{\rho} = \mathcal{T}_{\delta\bm{\rho}}(\delta\bm{\alpha},\delta\bm{\delta},\delta\dot{\bm{d}})\circ
\begin{Bmatrix}
\delta\bm{\alpha}\\
\delta\bm{\delta}\\
\left[\hat{\dot{d}}_1\right]-\hat{\dot{d}}_1
\vspace{0.2cm}\\
\left[\hat{\dot{d}}_N\right]-\hat{\dot{d}}_N
\end{Bmatrix}
\end{equation} 
\item[\label=\color{RoyalBlue} 8.] Compose $\left[\hat{\bm{r}}_1\right]$ and $\left[\hat{\bm{v}}_{1,\mathcal{L}}\right]$ with  Eq.~\eqref{eq:map}, i.e.
\begin{equation}
\begin{gathered}
\left[\hat{\bm{r}}_1\right] = \mathcal{T}_{\hat{\bm{r}}_1}\left(\delta\bm{\alpha},\delta\bm{\delta},\delta\bm{\rho}\right)\circ \delta\bm{\rho}\left(\delta\bm{\alpha},\delta\bm{\delta},\delta\dot{\bm{d}}\right) =
\mathcal{T}_{\hat{\bm{r}}_1}\left(\delta\bm{\alpha},\delta\bm{\delta},\delta\dot{\bm{d}}\right)\\
\left[\hat{\bm{v}}_{1,\mathcal{L}}\right] = \mathcal{T}_{\hat{\bm{v}}_{1,\mathcal{L}}}\left(\delta\bm{\alpha},\delta\bm{\delta},\delta\bm{\rho}\right)\circ \delta\bm{\rho}\left(\delta\bm{\alpha},\delta\bm{\delta},\delta\dot{\bm{d}}\right) =
\mathcal{T}_{\hat{\bm{v}}_{1,\mathcal{L}}}\left(\delta\bm{\alpha},\delta\bm{\delta},\delta\dot{\bm{d}}\right)
\end{gathered}
\end{equation}
\end{enumerate}

The described procedure provides an estimate for the object state $\hat{\bm{x}}_1$ at $t_1$. The state can be then converted to any desired representation $\hat{\bm{p}}_1$. In the adopted work, a reduced set of Keplerian parameters, namely semi-major axis $a$, eccentricity $e$, inclination $i$, right ascension of the ascending node $\Omega$ and argument of latitude $u$ is considered. Steps 3-8 of phase 3 are embedded into an \ac{ads} solver that automatically controls the accuracy of the resulting polynomial expansions and splits the uncertainty set, if required. The splitting decision is governed by tolerances set on the components of $\hat{\bm{p}}_1$~\citep{Wittig2015}. As a result, the final estimated state can be expressed as a manifold of Taylor expansions, i.e.
\begin{equation}
\label{eq:DAIOD_result}
[\hat{\bm{p}}_1] = \bigcup\limits_{i=1}^{N_s}\mathcal{T}_{\hat{\bm{p}}_1}\left(\delta\alpha_1,\delta\delta_1,\delta\dot{d}_1,\delta\alpha_N,\delta\delta_N,\delta\dot{d}_N\right)
\end{equation}

The state manifold given in Eq.~\eqref{eq:DAIOD_result} is an analytical mapping between the Doppler radar measurements space to the orbital state. This mapping is extremely useful for cataloguing purposes. For example, the solution corresponding to a generic set of deviations $\left\{\delta\alpha_1^*,\delta\delta_1^*,\delta\dot{d}_1^*,\delta\alpha_N^*,\delta\delta_N^*,\delta\dot{d}_N^*\right\}$ can be immediately obtained by the evaluation of  Eq.~\eqref{eq:DAIOD_result}, i.e., by
\begin{equation}
\hat{\bm{p}}_1^*=\bigcup\limits_{i=1}^{N_s}\mathcal{T}_{\hat{\bm{p}}_1}\left(\delta\alpha_1^*,\delta\delta_1^*,\delta\dot{d}_1^*,\delta\alpha_N^*,\delta\delta_N^*,\delta\dot{d}_N^*\right)
\end{equation}
The advantages granted by the polynomial description provided by Eq.~\eqref{eq:DAIOD_result} are not limited to polynomial evaluations. By applying the polynomial bounding techniques described in Section~\ref{sec:DA} one can easily estimate the bound of each component of the state vector $\hat{\bm{p}}_1$. Thus, an estimate of the envelope of the \textit{orbit set} compatible with the considered measurement uncertainties is readily  available.


\section{Performance assessment}
\label{sec:performance}
The performance of the DAIOD method is now investigated. The results of numerical simulations are illustrated in Section~\ref{subsec:Num_sim}, whereas Section~\ref{subsec:Real} is dedicated to the analysis of real data. For all DAIOD simulations, an expansion order of 4 and a maximum of 5 splits per direction were used. The tolerances for $a$, $e$, $i$, $\Omega$ and $u$ were set to $\{$\SI{0.01}{\kilo\meter}, 0.01, 1e-5~deg, 1e-5~deg, 1e-5~deg$\}$, respectively. All the simulations were run on an Intel i7-8565U CPU @ 1.80GHz and 16GB of RAM. The \ac{da} used in this work is CNES' PACE library.

\subsection{Numerical simulations}
\label{subsec:Num_sim}
\begin{figure*}[!t]
\centering
\subfloat[\label{subfig:Err_a_raw_5s}]{\includegraphics[trim=0cm 0cm 0cm 0.5cm, clip=true, width=0.5\textwidth]{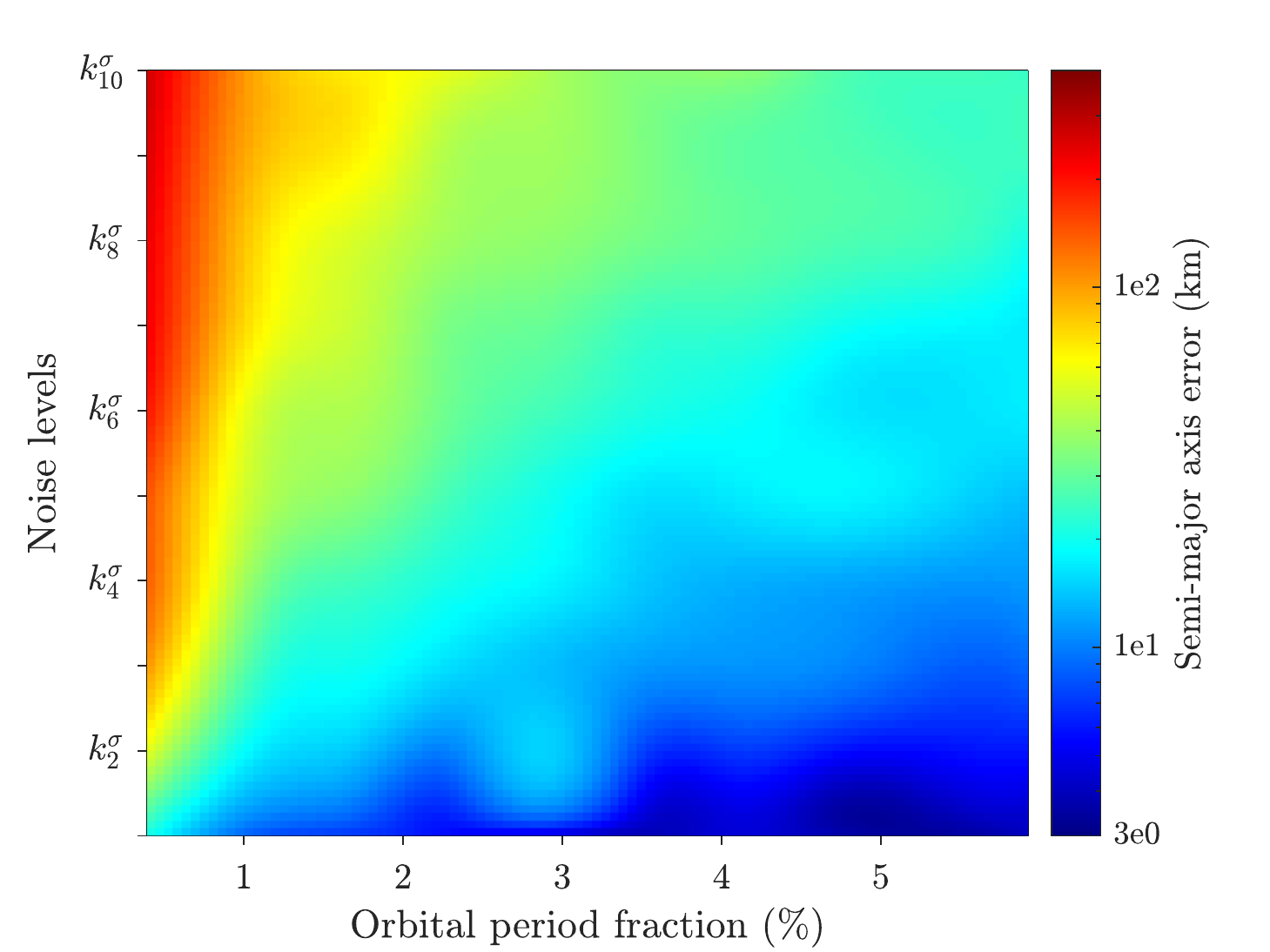}}
\subfloat[\label{subfig:B_a_raw_5s}]{\includegraphics[trim=0cm 0cm 0cm 0.5cm, clip=true, width=0.5\textwidth]{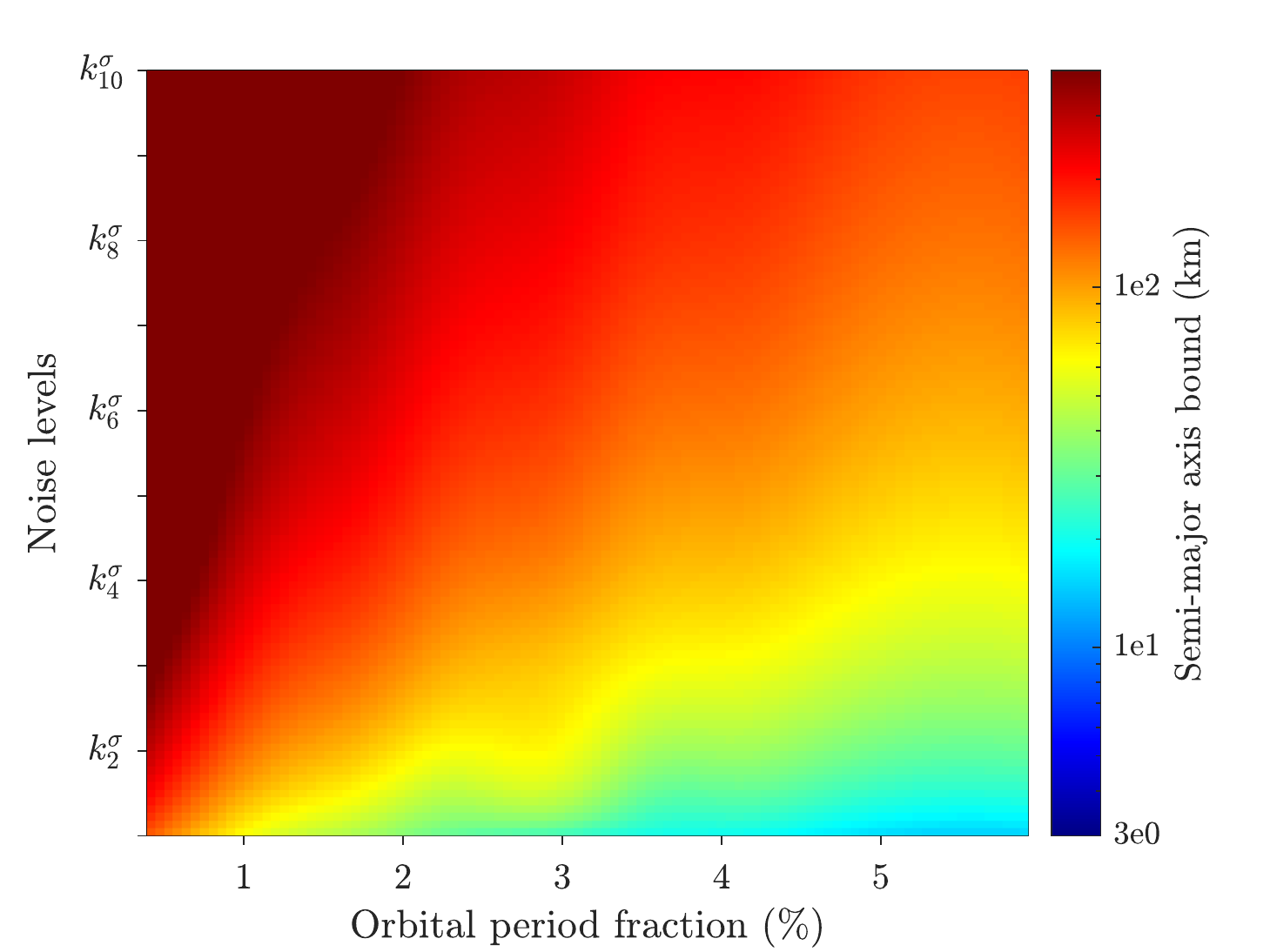}}\\
\caption{DAIOD performance in terms of mean error $\varepsilon_a$ \protect\subref{subfig:Err_a_raw_5s} and bound $b_a$ \protect\subref{subfig:B_a_raw_5s} in semi-major axis $a$ as a function of the fraction of observed orbital period and noise level $k^{\sigma}$ (raw data, $t_s$ = 5~s).}
\label{fig:Err_b_raw_5s}
\end{figure*}
The DAIOD method is first tested on numerical simulations targeting a subset of the NORAD \ac{leo} population. The simulations were carried out by downloading the latest \ac{tle} from the Space-Track\texorpdfstring{\protect\footnote{\href{https://www.space-track.org/}{https://www.space-track.org/}}}{} website, propagating them for one day with a high-fidelity dynamics (including the zonal/tesseral effects of the Earth's gravitational field up to order and degree 8, the gravitational pull of Sun and Moon, the solar radiation pressure and the drag perturbation), and considering as observer in surveillance mode the French bistatic radar sensor GRAVES~\citep{Muller2017}. An overall number of 2,000 passages was generated, and for each passage, 10 different measurement noise levels were considered, ranging from $k^{\sigma}_1$ = (10~mdeg; 10~mdeg; \SI{0.1}{\meter \per \second}) to $k^{\sigma}_{10}$ = (100~mdeg; 100~mdeg; \SI{1.0}{\meter \per \second}), with a step of (10~mdeg; 10~mdeg; \SI{0.1}{\meter \per \second}), where each triple indicates the noise standard deviations in azimuth, elevation, and range rate, respectively. The DAIOD method was then run on all the resulting 20,000 scenarios, and a comparison between estimated and real orbital elements performed.

\begin{table*}[!t]
    \centering
    \small
    \sisetup{round-mode=places,round-precision=4}
    \begin{tabular}{c c c c c c c}
    \hline\hline
    & & {$\Delta t_{obs}<0.01T$} & {$\Delta t_{obs}<0.02T$} & {$\Delta t_{obs}<0.03T$} & {$\Delta t_{obs}<0.04T$} & {$\Delta t_{obs}<0.05T$} \\ 
    \cline{2-7}
    		\multirow{5}{*}{\shortstack[c]{$\left[\varepsilon^{75\%}_a;b^{75\%}_a\right]$\\ (km)}}&	$k^{\sigma}_2$&	[\num[round-mode=figures,round-precision=5]{58.7292815607}\,;\,\num[round-mode=figures,round-precision=5]{401.9604606714}]&	[\num[round-mode=figures,round-precision=5]{39.4003737236}\,;\,\num[round-mode=figures,round-precision=5]{246.2131337377}]&	[\num[round-mode=figures,round-precision=5]{28.9338552028}\,;\,\num[round-mode=figures,round-precision=5]{189.2236380429}]&	[\num[round-mode=figures,round-precision=5]{22.9498883457}\,;\,\num[round-mode=figures,round-precision=5]{137.0833989709}]&	[\num[round-mode=figures,round-precision=5]{20.5698112765}\,;\,\num[round-mode=figures,round-precision=5]{124.4989405356}]\\
&	$k^{\sigma}_4$&	[\num[round-mode=figures,round-precision=5]{138.6539426122}\,;\,\num[round-mode=figures,round-precision=5]{860.2543636249}]&	[\num[round-mode=figures,round-precision=5]{77.0191828505}\,;\,\num[round-mode=figures,round-precision=5]{509.3951573384}]&	[\num[round-mode=figures,round-precision=5]{55.0917129594}\,;\,\num[round-mode=figures,round-precision=5]{394.6062869454}]&	[\num[round-mode=figures,round-precision=5]{43.5628003610}\,;\,\num[round-mode=figures,round-precision=5]{281.7427244721}]&	[\num[round-mode=figures,round-precision=5]{38.3440465049}\,;\,\num[round-mode=figures,round-precision=5]{253.3434565801}]\\
&	$k^{\sigma}_6$&	[\num[round-mode=figures,round-precision=5]{185.7844959520}\,;\,\num[round-mode=figures,round-precision=5]{1419.1136692140}]&	[\num[round-mode=figures,round-precision=5]{113.7332882354}\,;\,\num[round-mode=figures,round-precision=5]{786.3276140111}]&	[\num[round-mode=figures,round-precision=5]{86.7070650975}\,;\,\num[round-mode=figures,round-precision=5]{603.9365609954}]&	[\num[round-mode=figures,round-precision=5]{67.7319973685}\,;\,\num[round-mode=figures,round-precision=5]{432.0525453341}]&	[\num[round-mode=figures,round-precision=5]{60.4901959703}\,;\,\num[round-mode=figures,round-precision=5]{385.6796254854}]\\
&	$k^{\sigma}_8$&	[\num[round-mode=figures,round-precision=5]{251.9495966829}\,;\,\num[round-mode=figures,round-precision=5]{2034.6011193003}]&	[\num[round-mode=figures,round-precision=5]{166.6809129802}\,;\,\num[round-mode=figures,round-precision=5]{1055.4054850290}]&	[\num[round-mode=figures,round-precision=5]{116.4229872417}\,;\,\num[round-mode=figures,round-precision=5]{833.7985050244}]&	[\num[round-mode=figures,round-precision=5]{90.0794093864}\,;\,\num[round-mode=figures,round-precision=5]{588.7793097047}]&	[\num[round-mode=figures,round-precision=5]{77.6735134290}\,;\,\num[round-mode=figures,round-precision=5]{524.9536199722}]\\
&	$k^{\sigma}_{10}$&	[\num[round-mode=figures,round-precision=5]{281.8492814703}\,;\,\num[round-mode=figures,round-precision=5]{2869.1439639281}]&	[\num[round-mode=figures,round-precision=5]{189.6260585249}\,;\,\num[round-mode=figures,round-precision=5]{1364.8365969552}]&	[\num[round-mode=figures,round-precision=5]{141.0043054458}\,;\,\num[round-mode=figures,round-precision=5]{1060.3137886245}]&	[\num[round-mode=figures,round-precision=5]{105.9366592819}\,;\,\num[round-mode=figures,round-precision=5]{762.5941541499}]&	[\num[round-mode=figures,round-precision=5]{93.4113741166}\,;\,\num[round-mode=figures,round-precision=5]{670.8438283432}]\\
\hline
		\multirow{5}{*}{\shortstack[c]{$\left[\varepsilon^{75\%}_e;b^{75\%}_e\right]$\\ (-)}}&	$k^{\sigma}_2$&	[\num{0.0047782860}\,;\,\num{0.1632042581}]&	[\num{0.0032803612}\,;\,\num{0.1147652639}]&	[\num{0.0023787776}\,;\,\num{0.0870658782}]&	[\num{0.0018274855}\,;\,\num{0.0623272426}]&	[\num{0.0015798045}\,;\,\num{0.0503446557}]\\
&	$k^{\sigma}_4$&	[\num{0.0115644566}\,;\,\num{0.2542407480}]&	[\num{0.0071464086}\,;\,\num{0.2002452811}]&	[\num{0.0049113138}\,;\,\num{0.1607371732}]&	[\num{0.0038005362}\,;\,\num{0.1316727636}]&	[\num{0.0031675279}\,;\,\num{0.1158589515}]\\
&	$k^{\sigma}_6$&	[\num{0.0174190496}\,;\,\num{0.3283146345}]&	[\num{0.0104324226}\,;\,\num{0.2473473713}]&	[\num{0.0077742180}\,;\,\num{0.2086039329}]&	[\num{0.0061290766}\,;\,\num{0.1745721737}]&	[\num{0.0053046166}\,;\,\num{0.1576835332}]\\
&	$k^{\sigma}_8$&	[\num{0.0257964501}\,;\,\num{0.3901550838}]&	[\num{0.0153225961}\,;\,\num{0.2997791239}]&	[\num{0.0107125014}\,;\,\num{0.2490450958}]&	[\num{0.0082310368}\,;\,\num{0.2061193091}]&	[\num{0.0069642694}\,;\,\num{0.1922725717}]\\
&	$k^{\sigma}_{10}$&	[\num{0.0281883196}\,;\,\num{0.4571698525}]&	[\num{0.0180200680}\,;\,\num{0.3317263671}]&	[\num{0.0128938989}\,;\,\num{0.2843919533}]&	[\num{0.0095202543}\,;\,\num{0.2394830017}]&	[\num{0.0084590727}\,;\,\num{0.2184525527}]\\
\hline
		\multirow{5}{*}{\shortstack[c]{$\left[\varepsilon^{75\%}_i;b^{75\%}_i\right]$\\ (deg)}}&	$k^{\sigma}_2$&	[\num{0.1518449443}\,;\,\num{0.7576721098}]&	[\num{0.0995593046}\,;\,\num{0.4687036852}]&	[\num{0.0698645778}\,;\,\num{0.3679045598}]&	[\num{0.0490416498}\,;\,\num{0.2553859252}]&	[\num{0.0442297958}\,;\,\num{0.2339811366}]\\
&	$k^{\sigma}_4$&	[\num{0.2658264533}\,;\,\num{1.5171114406}]&	[\num{0.1822556382}\,;\,\num{0.9261657790}]&	[\num{0.1360656887}\,;\,\num{0.7436712170}]&	[\num{0.0972543677}\,;\,\num{0.5150166850}]&	[\num{0.0824909076}\,;\,\num{0.4711196576}]\\
&	$k^{\sigma}_6$&	[\num{0.3658361087}\,;\,\num{2.3010754004}]&	[\num{0.2516655029}\,;\,\num{1.3406715365}]&	[\num{0.1893061017}\,;\,\num{1.1134532742}]&	[\num{0.1444774597}\,;\,\num{0.7700992904}]&	[\num{0.1247910624}\,;\,\num{0.7018157181}]\\
&	$k^{\sigma}_8$&	[\num{0.6239324407}\,;\,\num{2.9912066237}]&	[\num{0.3827451063}\,;\,\num{1.7692313903}]&	[\num{0.2658730218}\,;\,\num{1.4749547229}]&	[\num{0.1942954454}\,;\,\num{1.0326639424}]&	[\num{0.1736382711}\,;\,\num{0.9302831962}]\\
&	$k^{\sigma}_{10}$&	[\num{0.7136951352}\,;\,\num{3.7427067073}]&	[\num{0.4366552355}\,;\,\num{2.1525201264}]&	[\num{0.3131736349}\,;\,\num{1.8211547167}]&	[\num{0.2332137865}\,;\,\num{1.3000764011}]&	[\num{0.2065231459}\,;\,\num{1.1584411512}]\\
\hline
		\multirow{5}{*}{\shortstack[c]{$\left[\varepsilon^{75\%}_{\Omega};b^{75\%}_{\Omega}\right]$\\ (deg)}}&	$k^{\sigma}_2$&	[\num{0.1593780015}\,;\,\num{0.5821815201}]&	[\num{0.0992962637}\,;\,\num{0.5339264953}]&	[\num{0.0754936781}\,;\,\num{0.3796535875}]&	[\num{0.0585246822}\,;\,\num{0.2759481247}]&	[\num{0.0511930714}\,;\,\num{0.2601234567}]\\
&	$k^{\sigma}_4$&	[\num{0.3079846662}\,;\,\num{1.1969361663}]&	[\num{0.1988772776}\,;\,\num{1.0837265319}]&	[\num{0.1453302106}\,;\,\num{0.7539248218}]&	[\num{0.1095144566}\,;\,\num{0.5554617344}]&	[\num{0.0950042513}\,;\,\num{0.5220203781}]\\
&	$k^{\sigma}_6$&	[\num{0.4016891315}\,;\,\num{1.8027173862}]&	[\num{0.2860681816}\,;\,\num{1.6403480747}]&	[\num{0.2253308136}\,;\,\num{1.1006709803}]&	[\num{0.1755921250}\,;\,\num{0.8360760870}]&	[\num{0.1477488309}\,;\,\num{0.7848409017}]\\
&	$k^{\sigma}_8$&	[\num{0.5904108656}\,;\,\num{2.4218894965}]&	[\num{0.3910107713}\,;\,\num{2.1696947349}]&	[\num{0.2889718543}\,;\,\num{1.4742882519}]&	[\num{0.2225122145}\,;\,\num{1.1186865374}]&	[\num{0.1932639475}\,;\,\num{1.0466048371}]\\
&	$k^{\sigma}_{10}$&	[\num{0.6771224771}\,;\,\num{3.0045421683}]&	[\num{0.4464837942}\,;\,\num{2.7049871338}]&	[\num{0.3537266480}\,;\,\num{1.8130315498}]&	[\num{0.2706035075}\,;\,\num{1.4068862634}]&	[\num{0.2399929710}\,;\,\num{1.3065396130}]\\
\hline
		\multirow{5}{*}{\shortstack[c]{$\left[\varepsilon^{75\%}_u;b^{75\%}_u\right]$\\ (deg)}}&	$k^{\sigma}_2$&	[\num{0.0718678539}\,;\,\num{0.3510024130}]&	[\num{0.0498490754}\,;\,\num{0.2197903433}]&	[\num{0.0357409687}\,;\,\num{0.1971861863}]&	[\num{0.0264321724}\,;\,\num{0.1489706356}]&	[\num{0.0235334448}\,;\,\num{0.1253365728}]\\
&	$k^{\sigma}_4$&	[\num{0.1483047505}\,;\,\num{0.7076304959}]&	[\num{0.0924373973}\,;\,\num{0.4525586937}]&	[\num{0.0704685577}\,;\,\num{0.4014694726}]&	[\num{0.0549660734}\,;\,\num{0.3014278939}]&	[\num{0.0471331044}\,;\,\num{0.2526993428}]\\
&	$k^{\sigma}_6$&	[\num{0.2129294508}\,;\,\num{1.0547541711}]&	[\num{0.1408311072}\,;\,\num{0.6954003052}]&	[\num{0.1083820864}\,;\,\num{0.6163429683}]&	[\num{0.0806708998}\,;\,\num{0.4542487304}]&	[\num{0.0710519293}\,;\,\num{0.3775360737}]\\
&	$k^{\sigma}_8$&	[\num{0.2853383596}\,;\,\num{1.4269319701}]&	[\num{0.1821557785}\,;\,\num{0.9338298977}]&	[\num{0.1396378177}\,;\,\num{0.8351098431}]&	[\num{0.1044930520}\,;\,\num{0.6066646864}]&	[\num{0.0940215938}\,;\,\num{0.5038818247}]\\
&	$k^{\sigma}_{10}$&	[\num{0.3596983111}\,;\,\num{1.8105414969}]&	[\num{0.2351070167}\,;\,\num{1.2167365115}]&	[\num{0.1753958075}\,;\,\num{1.0410977219}]&	[\num{0.1324687717}\,;\,\num{0.7614602033}]&	[\num{0.1147050157}\,;\,\num{0.6329369681}]\\
    \hline\hline
    \end{tabular}
    \vspace{0.2cm}
    \caption{DAIOD performance expressed in terms of 75th percentiles of the errors $\varepsilon$ and bounds $b$ in semi-major axis $a$, eccentricity $e$, inclination $i$, right ascension of the ascending node $\Omega$ and argument of latitude $u$ as a function of the noise level $k^{\sigma}$ and the fraction of observed orbital period $\Delta t_{obs}$ (raw data, $t_s$ = 5~s).}
    \label{tab:Results_err_b_raw_5s}
\end{table*}

Figure~\ref{fig:Err_b_raw_5s} shows the performance of the DAIOD method in terms of mean errors (Fig.~\ref{subfig:Err_a_raw_5s}) and bounds (Fig.~\ref{subfig:B_a_raw_5s}) of the semi-major axis estimate when processing raw data with a sampling time $t_s$ of 5~s. The analysis is performed against noise levels and observation length by assigning each single passage to the corresponding cell, and then computing the average of the resulting errors and bounds per cell. Let us first analyse Fig.~\ref{subfig:Err_a_raw_5s}. The plot shows the trend of the semi-major axis error $\varepsilon_a$ as a function of the noise level $k^{\sigma}$ and the fraction of orbital period observed. Two major trends can be identified. If one considers a specific noise level, the accuracy of the obtained estimate increases for increasing portions of observed arc. This is expected, as a better estimation of the orbit curvature can be obtained when processing a longer arc. On the other hand, an increase in the estimate error for increasing noise levels is noticed. Overall, three main regions can be identified. The first region covers the bottom-right portion of the plot, and collects longer passages with relatively low noise levels. Here $\varepsilon_a$ is in average around few kilometers. The second region is located in the top-left corner, and it is characterised by short passages with increasingly higher noise levels. The combination of very short-arc observations with large noise values is the most critical condition for the DAIOD method, with semi-major axis errors ranging from \SI{80}{\kilo \meter} to about \SI{400}{\kilo \meter}. The third region, instead, covers the remaining portion of the map, with $\varepsilon_a \in [10, 80]$\ km. 

Figure~\ref{subfig:B_a_raw_5s} shows the trend of the estimated semi-major axis bounds $b_a$. The trend closely resembles the one of Fig.~\ref{subfig:Err_a_raw_5s}. Furthermore, the comparison between values of $\varepsilon_a$ and $b_a$ for a generic noise-arc length couple clearly shows that the estimated bound $b_a$ is systematically larger than the resulting error $\varepsilon_a$. That is, an estimate may be inaccurate, but the uncertainty associated with it accurately reproduces this inaccuracy. This is an important result, as it proves that the DAIOD method provides reliable estimates regardless of the observation conditions.

\begin{table*}[!t]
    \centering
    \small
    \sisetup{round-mode=places,round-precision=4}
    \begin{tabular}{l c c c c c}
    \hline\hline
    & {$\Delta t_{obs}<0.01T$} & {$\Delta t_{obs}<0.02T$} & {$\Delta t_{obs}<0.03T$} & {$\Delta t_{obs}<0.04T$} & {$\Delta t_{obs}<0.05T$} \\ 
    \cline{2-6}
		Centre-only&	\num{0.439882697947214}&	\num{0.555393586005831}&	\num{0.669201520912548}&	\num{0.762790697674418}&	\num{0.800334634690463}\\
		Centre+64 corners&	\num{0.967741935483871}&	\num{0.983965014577260}&	\num{0.989543726235742}&	\num{0.991362126245847}&	\num{0.992749581706637}\\
    \hline\hline
    \end{tabular}
    \vspace{0.2cm}
    \caption{DAIOD success rate as function of the fraction of observed orbital period $\Delta t_{obs}$ when computing Gauss' guess with the centre of the angles uncertainty sets only (line 1), or by adding the 64 corners (line 2)(raw data, $t_s$ = 5~s, $k^{\sigma}_{10}$).
    \label{tab:Results_success_raw_5s}}
\end{table*}

Table~\ref{tab:Results_err_b_raw_5s} offers a complete picture of the DAIOD performance when processing raw data with a sampling time of 5~s. The table shows the $Q_3$ quartiles (75th percentiles) of the errors and bounds in $a$, $e$, $i$, $\Omega$ and $u$ as a function of the noise level $k^{\sigma}$ and the observation interval $\Delta t_{obs}$, as expressed in terms of fraction of the orbital period $T$. As expected, the DAIOD estimate accuracy increases for increasing arc lengths and decreasing noise levels, while the estimated bounds $b$ are always larger than the actual errors $\varepsilon$. Overall, a general reduction of one order of magnitude in both $\varepsilon$ and $b$ can be noticed while passing from the top-left to the bottom-right corner of the noise-arc length plane. 

While the performance of the DAIOD method is strongly dependent on the duration of the passage and the noise level, its success rate, instead, is always very high and barely affected by the observation conditions. This second important result is illustrated by the second row of Table~\ref{tab:Results_success_raw_5s}. The table shows the trend of the DAIOD success rate as a function of the arc length for the highest considered noise level ($k^{\sigma}_{10}$). The success rate remains almost unaltered while moving from the top-left to the top-right of the noise-arc length plane, passing from 0.9677 to 0.9927. This robustness is enabled by the scanning of the angular uncertainty set performed at the beginning of phase 2, as confirmed by the data in the first row. The line reports the success rate when only the nominal angular measurements (centre-only) are used to solve Gauss' problem. The improvement granted by the adopted approach is remarkable, especially when the passage is very short. This last aspect, i.e. the capability of handling very short-arc observations and providing an estimate regardless of the noise level represents the third important result for the method.

\begin{figure*}[!t]
\centering
\subfloat[\label{subfig:Sets_raw_5s}]{\includegraphics[trim=0cm 0cm 0cm 0cm, clip=true, width=0.5\textwidth]{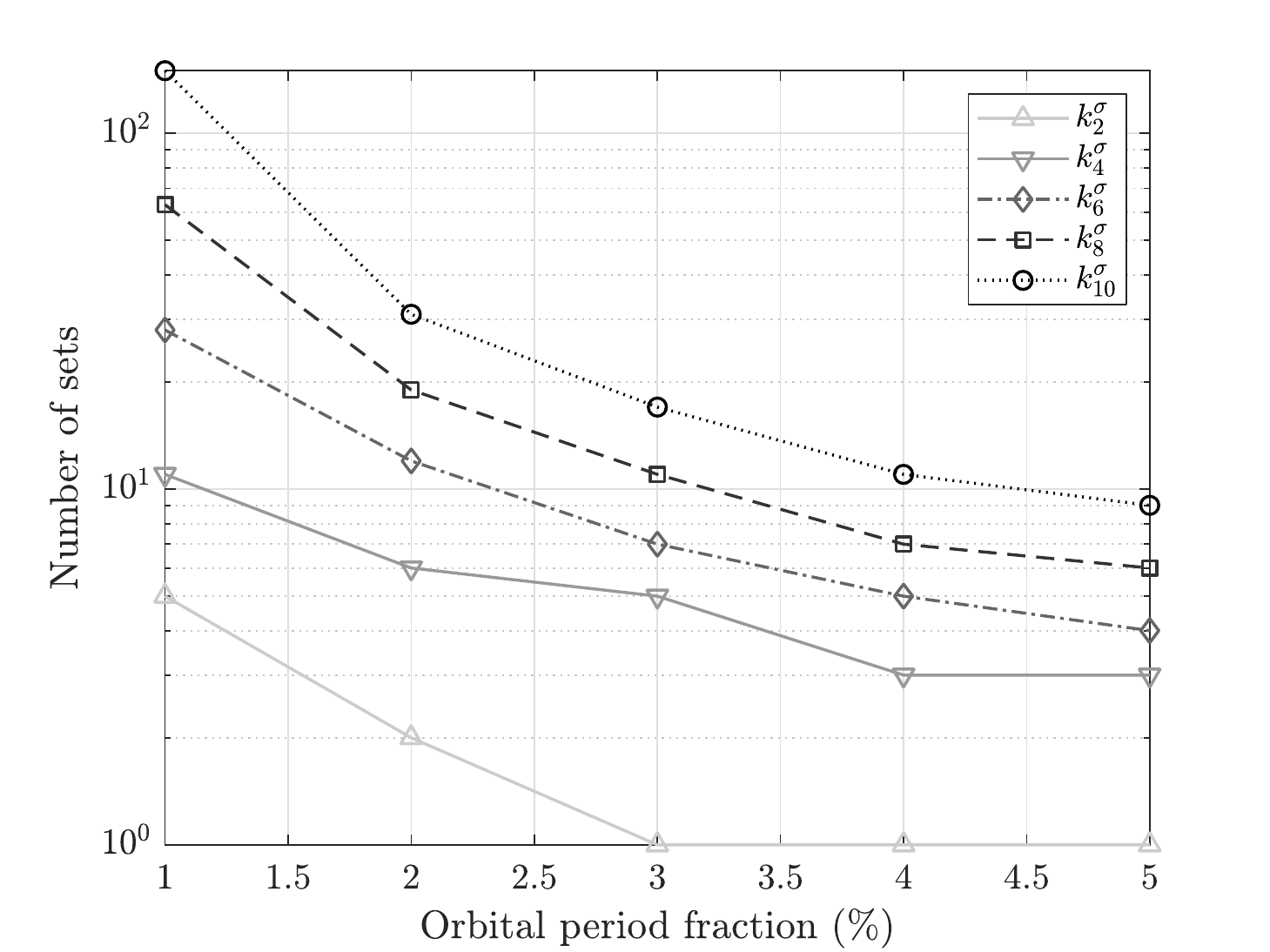}}
\subfloat[\label{subfig:Time_raw_5s}]{\includegraphics[trim=0cm 0cm 0cm 0cm, clip=true, width=0.5\textwidth]{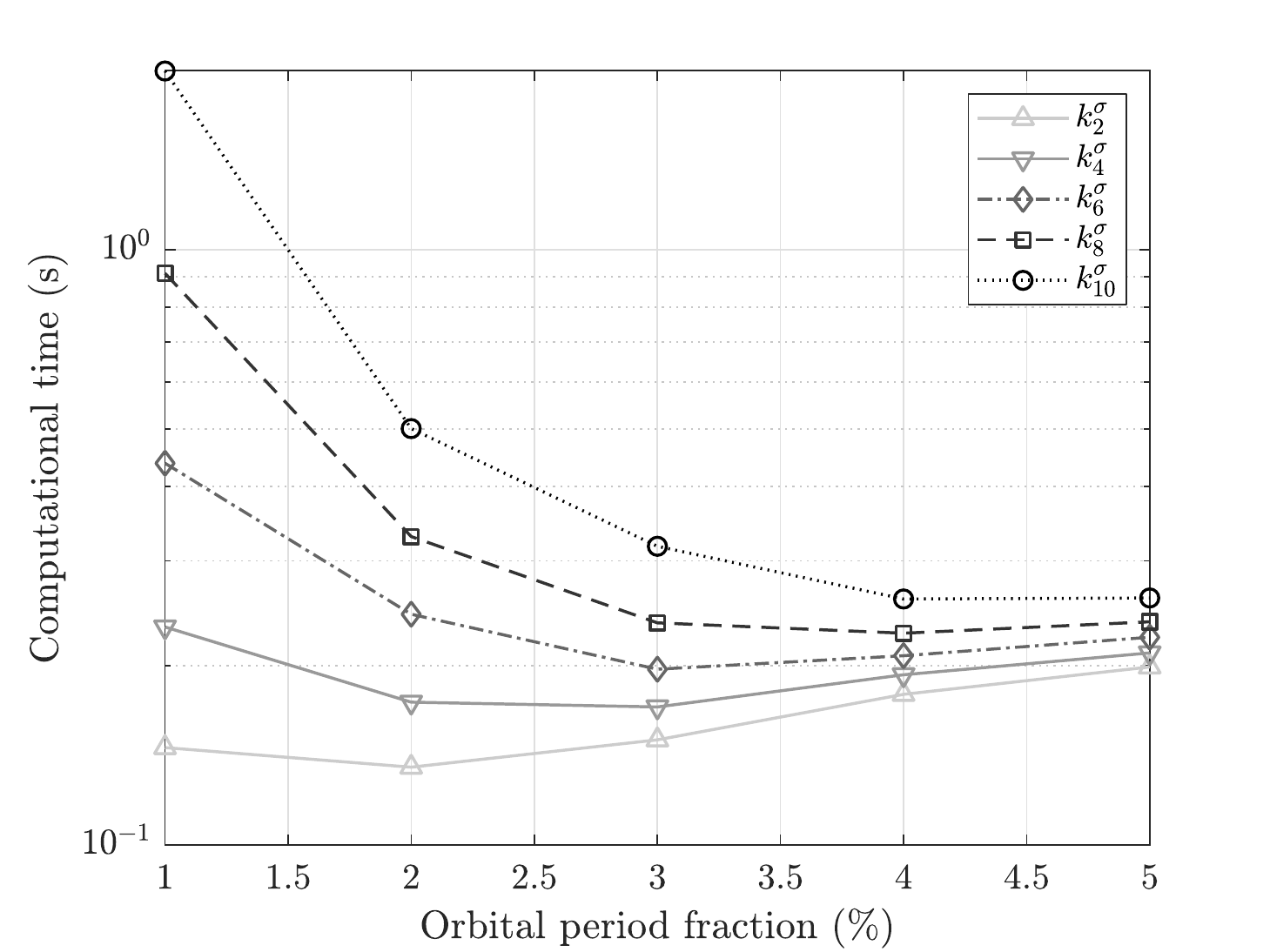}}\\
\caption{DAIOD performance in terms of 75th percentile of the required number of sets \protect\subref{subfig:Sets_raw_5s} and computational time \protect\subref{subfig:Time_raw_5s}  as a function of the fraction of observed orbital period and noise level $k^{\sigma}$ (raw data, $t_s$ = 5~s).}
\label{fig:Time_sets_raw_5s}
\end{figure*}

\begin{table*}[!t]
    \centering
    \small
    \sisetup{round-mode=places,round-precision=4}
    \begin{tabular}{c c c c c c c}
    \hline\hline
    & & {$\Delta t_{obs}<0.01T$} & {$\Delta t_{obs}<0.02T$} & {$\Delta t_{obs}<0.03T$} & {$\Delta t_{obs}<0.04T$} & {$\Delta t_{obs}<0.05T$} \\ 
    \cline{2-7}
    		\multirow{5}{*}{\shortstack[c]{$\left[\varepsilon^{75\%}_a;b^{75\%}_a\right]$\\ (km)\\$t_s$ = 3~s}}&	$k^{\sigma}_2$&	[\num[round-mode=figures,round-precision=5]{43.3802729245}\,;\,\num[round-mode=figures,round-precision=5]{276.1135902785}]&	[\num[round-mode=figures,round-precision=5]{29.4110369847}\,;\,\num[round-mode=figures,round-precision=5]{166.6813948364}]&	[\num[round-mode=figures,round-precision=5]{23.7680143432}\,;\,\num[round-mode=figures,round-precision=5]{122.9119201515}]&	[\num[round-mode=figures,round-precision=5]{18.8030675414}\,;\,\num[round-mode=figures,round-precision=5]{83.5440720981}]&	[\num[round-mode=figures,round-precision=5]{16.4329300930}\,;\,\num[round-mode=figures,round-precision=5]{73.5645087147}]\\
&	$k^{\sigma}_4$&	[\num[round-mode=figures,round-precision=5]{70.6484864544}\,;\,\num[round-mode=figures,round-precision=5]{559.0625244999}]&	[\num[round-mode=figures,round-precision=5]{47.6207182308}\,;\,\num[round-mode=figures,round-precision=5]{295.1085491374}]&	[\num[round-mode=figures,round-precision=5]{41.3700494166}\,;\,\num[round-mode=figures,round-precision=5]{232.1114860416}]&	[\num[round-mode=figures,round-precision=5]{34.1340551744}\,;\,\num[round-mode=figures,round-precision=5]{163.9559659931}]&	[\num[round-mode=figures,round-precision=5]{30.7461587827}\,;\,\num[round-mode=figures,round-precision=5]{145.2366782547}]\\
&	$k^{\sigma}_6$&	[\num[round-mode=figures,round-precision=5]{109.7883155838}\,;\,\num[round-mode=figures,round-precision=5]{942.8077990331}]&	[\num[round-mode=figures,round-precision=5]{70.2173649452}\,;\,\num[round-mode=figures,round-precision=5]{457.7475687645}]&	[\num[round-mode=figures,round-precision=5]{60.3741517365}\,;\,\num[round-mode=figures,round-precision=5]{344.0798493090}]&	[\num[round-mode=figures,round-precision=5]{51.1874442429}\,;\,\num[round-mode=figures,round-precision=5]{249.8904711084}]&	[\num[round-mode=figures,round-precision=5]{45.3087764810}\,;\,\num[round-mode=figures,round-precision=5]{217.3791331032}]\\
&	$k^{\sigma}_8$&	[\num[round-mode=figures,round-precision=5]{137.2160235453}\,;\,\num[round-mode=figures,round-precision=5]{1173.1402499476}]&	[\num[round-mode=figures,round-precision=5]{89.7284800230}\,;\,\num[round-mode=figures,round-precision=5]{574.4439803421}]&	[\num[round-mode=figures,round-precision=5]{72.8729407425}\,;\,\num[round-mode=figures,round-precision=5]{445.5319278680}]&	[\num[round-mode=figures,round-precision=5]{60.2221800968}\,;\,\num[round-mode=figures,round-precision=5]{333.5798865427}]&	[\num[round-mode=figures,round-precision=5]{53.6548288413}\,;\,\num[round-mode=figures,round-precision=5]{295.9907518570}]\\
&	$k^{\sigma}_{10}$&	[\num[round-mode=figures,round-precision=5]{172.0953158922}\,;\,\num[round-mode=figures,round-precision=5]{1564.1993232217}]&	[\num[round-mode=figures,round-precision=5]{104.6004190845}\,;\,\num[round-mode=figures,round-precision=5]{730.4171761650}]&	[\num[round-mode=figures,round-precision=5]{91.8224209298}\,;\,\num[round-mode=figures,round-precision=5]{562.0778808556}]&	[\num[round-mode=figures,round-precision=5]{77.0113476380}\,;\,\num[round-mode=figures,round-precision=5]{431.3594236819}]&	[\num[round-mode=figures,round-precision=5]{69.8970952773}\,;\,\num[round-mode=figures,round-precision=5]{375.4608879351}]\\
\hline
			\multirow{5}{*}{\shortstack[c]{$\left[\varepsilon^{75\%}_a;b^{75\%}_a\right]$\\ (km)\\$t_s$ = 5~s}}&	$k^{\sigma}_2$&	[\num[round-mode=figures,round-precision=5]{48.0798589443}\,;\,\num[round-mode=figures,round-precision=5]{408.9782702055}]&	[\num[round-mode=figures,round-precision=5]{33.4736296940}\,;\,\num[round-mode=figures,round-precision=5]{222.5919620799}]&	[\num[round-mode=figures,round-precision=5]{24.2400143035}\,;\,\num[round-mode=figures,round-precision=5]{143.4924481430}]&	[\num[round-mode=figures,round-precision=5]{19.5314909143}\,;\,\num[round-mode=figures,round-precision=5]{98.0166025365}]&	[\num[round-mode=figures,round-precision=5]{17.5923463075}\,;\,\num[round-mode=figures,round-precision=5]{86.8528285309}]\\
&	$k^{\sigma}_4$&	[\num[round-mode=figures,round-precision=5]{95.9366468247}\,;\,\num[round-mode=figures,round-precision=5]{775.7642410008}]&	[\num[round-mode=figures,round-precision=5]{61.8545576313}\,;\,\num[round-mode=figures,round-precision=5]{399.0837959645}]&	[\num[round-mode=figures,round-precision=5]{47.3270834739}\,;\,\num[round-mode=figures,round-precision=5]{298.8228556327}]&	[\num[round-mode=figures,round-precision=5]{36.9259900972}\,;\,\num[round-mode=figures,round-precision=5]{196.2662699330}]&	[\num[round-mode=figures,round-precision=5]{32.8052157304}\,;\,\num[round-mode=figures,round-precision=5]{169.6474519599}]\\
&	$k^{\sigma}_6$&	[\num[round-mode=figures,round-precision=5]{126.4035721424}\,;\,\num[round-mode=figures,round-precision=5]{1130.9291791947}]&	[\num[round-mode=figures,round-precision=5]{86.8142306841}\,;\,\num[round-mode=figures,round-precision=5]{582.4671181268}]&	[\num[round-mode=figures,round-precision=5]{69.6351853172}\,;\,\num[round-mode=figures,round-precision=5]{435.2944120086}]&	[\num[round-mode=figures,round-precision=5]{55.7689492555}\,;\,\num[round-mode=figures,round-precision=5]{307.7469004554}]&	[\num[round-mode=figures,round-precision=5]{48.8570837077}\,;\,\num[round-mode=figures,round-precision=5]{257.5030873290}]\\
&	$k^{\sigma}_8$&	[\num[round-mode=figures,round-precision=5]{204.2552963551}\,;\,\num[round-mode=figures,round-precision=5]{1600.4794303829}]&	[\num[round-mode=figures,round-precision=5]{113.8557055605}\,;\,\num[round-mode=figures,round-precision=5]{772.5094795756}]&	[\num[round-mode=figures,round-precision=5]{90.0569523674}\,;\,\num[round-mode=figures,round-precision=5]{582.4788869356}]&	[\num[round-mode=figures,round-precision=5]{72.0548540043}\,;\,\num[round-mode=figures,round-precision=5]{423.1913718525}]&	[\num[round-mode=figures,round-precision=5]{63.3017357492}\,;\,\num[round-mode=figures,round-precision=5]{344.7606052575}]\\
&	$k^{\sigma}_{10}$&	[\num[round-mode=figures,round-precision=5]{224.5260813060}\,;\,\num[round-mode=figures,round-precision=5]{2201.3604384754}]&	[\num[round-mode=figures,round-precision=5]{129.8877501286}\,;\,\num[round-mode=figures,round-precision=5]{961.4785919084}]&	[\num[round-mode=figures,round-precision=5]{109.5498581762}\,;\,\num[round-mode=figures,round-precision=5]{720.1589047938}]&	[\num[round-mode=figures,round-precision=5]{85.4568109685}\,;\,\num[round-mode=figures,round-precision=5]{517.3218951022}]&	[\num[round-mode=figures,round-precision=5]{77.5914541330}\,;\,\num[round-mode=figures,round-precision=5]{440.3750425059}]\\
\hline
			\multirow{5}{*}{\shortstack[c]{$\left[\varepsilon^{75\%}_a;b^{75\%}_a\right]$\\ (km)\\$t_s$ = 7.5~s}}&	$k^{\sigma}_2$&	[\num[round-mode=figures,round-precision=5]{46.4813527830}\,;\,\num[round-mode=figures,round-precision=5]{647.0663955394}]&	[\num[round-mode=figures,round-precision=5]{33.6351405446}\,;\,\num[round-mode=figures,round-precision=5]{309.9133944692}]&	[\num[round-mode=figures,round-precision=5]{27.4462961909}\,;\,\num[round-mode=figures,round-precision=5]{183.7859613451}]&	[\num[round-mode=figures,round-precision=5]{21.2516931425}\,;\,\num[round-mode=figures,round-precision=5]{118.8818878029}]&	[\num[round-mode=figures,round-precision=5]{18.1298709051}\,;\,\num[round-mode=figures,round-precision=5]{103.2891879067}]\\
&	$k^{\sigma}_4$&	[\num[round-mode=figures,round-precision=5]{104.9060893227}\,;\,\num[round-mode=figures,round-precision=5]{1256.7077427853}]&	[\num[round-mode=figures,round-precision=5]{67.3912211288}\,;\,\num[round-mode=figures,round-precision=5]{554.4730879543}]&	[\num[round-mode=figures,round-precision=5]{54.1781970611}\,;\,\num[round-mode=figures,round-precision=5]{385.2751453495}]&	[\num[round-mode=figures,round-precision=5]{43.0619442665}\,;\,\num[round-mode=figures,round-precision=5]{238.1090077361}]&	[\num[round-mode=figures,round-precision=5]{38.1542396581}\,;\,\num[round-mode=figures,round-precision=5]{191.4295527832}]\\
&	$k^{\sigma}_6$&	[\num[round-mode=figures,round-precision=5]{157.9093945418}\,;\,\num[round-mode=figures,round-precision=5]{1954.2998420584}]&	[\num[round-mode=figures,round-precision=5]{91.7478395604}\,;\,\num[round-mode=figures,round-precision=5]{780.1624566646}]&	[\num[round-mode=figures,round-precision=5]{71.8966958462}\,;\,\num[round-mode=figures,round-precision=5]{539.8253075326}]&	[\num[round-mode=figures,round-precision=5]{58.0664230485}\,;\,\num[round-mode=figures,round-precision=5]{364.7123103082}]&	[\num[round-mode=figures,round-precision=5]{52.5565660392}\,;\,\num[round-mode=figures,round-precision=5]{299.5486078888}]\\
&	$k^{\sigma}_8$&	[\num[round-mode=figures,round-precision=5]{192.1113300196}\,;\,\num[round-mode=figures,round-precision=5]{2900.9690774847}]&	[\num[round-mode=figures,round-precision=5]{120.6663057543}\,;\,\num[round-mode=figures,round-precision=5]{1004.1694994716}]&	[\num[round-mode=figures,round-precision=5]{95.1372908410}\,;\,\num[round-mode=figures,round-precision=5]{710.8998542922}]&	[\num[round-mode=figures,round-precision=5]{74.7566894762}\,;\,\num[round-mode=figures,round-precision=5]{488.5271097226}]&	[\num[round-mode=figures,round-precision=5]{66.1283572803}\,;\,\num[round-mode=figures,round-precision=5]{408.7458011760}]\\
&	$k^{\sigma}_{10}$&	[\num[round-mode=figures,round-precision=5]{233.8451005369}\,;\,\num[round-mode=figures,round-precision=5]{2895.5401148811}]&	[\num[round-mode=figures,round-precision=5]{138.6891925877}\,;\,\num[round-mode=figures,round-precision=5]{1261.9210155116}]&	[\num[round-mode=figures,round-precision=5]{114.0034116965}\,;\,\num[round-mode=figures,round-precision=5]{929.7147416377}]&	[\num[round-mode=figures,round-precision=5]{93.4731841652}\,;\,\num[round-mode=figures,round-precision=5]{603.9024687032}]&	[\num[round-mode=figures,round-precision=5]{80.4965289981}\,;\,\num[round-mode=figures,round-precision=5]{480.0587625922}]\\
    \hline\hline
    \end{tabular}
    \vspace{0.2cm}
    \caption{DAIOD performance expressed in terms of 75th percentiles of the errors $\varepsilon$ and bounds $b$ in semi-major axis $a$ as a function of the noise level $k^{\sigma}$ and the fraction of observed orbital period $\Delta t_{obs}$ when performing measurement regression with three sampling times: 3~s, 5~s, 7.5~s.}
    \label{tab:Results_regr}
\end{table*}

The DAIOD algorithm exploits the \ac{ads} technique to describe the solution uncertainty region accurately. Figure~\ref{subfig:Sets_raw_5s} shows the cardinality of the resulting ADS manifold as a function of the arc length and noise level $k^{\sigma}$. The number of sets is affected by both factors. More specifically, the value decreases for increasing arc lengths and decreasing noise levels, and can reach the ideal single-set description in case of low-noise measurements covering a relatively large portion of the object trajectory ($k^{\sigma}_2$, $\Delta t_{obs}\geq 0.03T$). This aspect has a direct consequence on the required computational time $t_{CPU}$ (see Fig.~\ref{subfig:Time_raw_5s}), for which a non-monotonic trend is identified. More specifically, for all the considered noise levels, $t_{CPU}$ first decreases and then increases for increasing values of $\Delta t_{obs}$.  This behaviour is determined by two factors: the measurement processing of phase 1 and the \ac{ads}. When the arc is short, the impact of phase 1 is negligible, so the overall computational time is mainly governed by phase 3, which generates \ac{ads} sets whose number decreases for increasing arc lengths. When the processed arc is longer, the number of sets is lower, thus the time required by the \ac{ads}  becomes comparable with the one dedicated to phase 1, which increases for longer arc lengths. Overall, fractions of seconds per passage are generally required.

\begin{table*}[!t]
    \centering
    \small
    \sisetup{round-mode=places,round-precision=4}
    \begin{tabular}{l c c c c c}
    \hline\hline
    & {$\Delta t_{obs}<0.01T$} & {$\Delta t_{obs}<0.02T$} & {$\Delta t_{obs}<0.03T$} & {$\Delta t_{obs}<0.04T$} & {$\Delta t_{obs}<0.05T$} \\ 
    \cline{2-6}
		$\varepsilon^{75\%}_a$ (km)& \num[round-mode=figures,round-precision=5]{4455.0377736078}&	\num[round-mode=figures,round-precision=5]{4178.3413000632}&	\num[round-mode=figures,round-precision=5]{1890.9406480951}&	\num[round-mode=figures,round-precision=5]{851.9506891696}&	\num[round-mode=figures,round-precision=5]{564.4639858670}\\
		$\varepsilon^{75\%}_e$ (-)&	\num{0.9355083771}&	\num{0.7984429760}&	\num{0.2147154021}&	\num{0.0866403695}&	\num{0.0596468234}\\
		$\varepsilon^{75\%}_i$ (deg)&		\num{21.1941522230}&	\num{3.5052160102}&	\num{1.0434308150}&	\num{0.4310718836}&	\num{0.2703258326}\\
		$\varepsilon^{75\%}_{\Omega}$ (deg)&	\num{41.9796054191}&	\num{5.2205092720}&	\num{0.8501498341}&	\num{0.3170243075}&	\num{0.2291048808}\\
		$\varepsilon^{75\%}_u$ (deg)&	\num{18.6161982045}&	\num{2.9499357914}&	\num{0.7055527423}&	\num{0.2422592587}&	\num{0.1590516901}\\
    \hline\hline
    \end{tabular}
    \vspace{0.2cm}
    \caption{\ac{dim} performance expressed in terms of 75th percentiles of the errors $\varepsilon$ in Keplerian parameters as a function of the fraction of observed orbital period $\Delta t_{obs}$ (regressed measurements, $t_s$ = 5~s, $k^{\sigma}_2$).}
    \label{tab:Results_DIM}
\end{table*}

All the analyses presented so far were done considering the DAIOD performance when processing raw data with a sampling time of 5 s. Table~\ref{tab:Results_regr} shows the results obtained when performing measurement regression with different values of $t_s$. More specifically, the table shows the 75th percentile of the errors $\varepsilon_a$ and bounds $b_a$ in semi-major axis as a function of noise level $k^{\sigma}$ and observed arc length $\Delta t_{obs}$ for regressed measurements obtained for three different values of $t_s$: 3~s, 5~s, and 7.5~s. Let us first analyse the effect of measurement regression, thus comparing the $t_s$ = 5~s block with its counterpart in Table~\ref{tab:Results_err_b_raw_5s}. The introduction of measurement regression grants a general reduction of both errors and bounds, which is evident for all the noise levels and the observed arc lengths. This is an expected outcome, as the regression allows the solver to mitigate the effect of measurement noise. Similarly, a decrease in the sampling time (i.e., an increase in the number of processed measurements) produces more accurate results and smaller bounds. As an example, if we consider the ($k^{\sigma}_{10}$; $\Delta t_{obs}<$ 0.01 $T$) critical case and compare the results with $t_{s}$ = 5~s and $t_{s}$ = 3~s, a drop in the semi-major axis error can be noticed, passing from 224~km to 172~km. The reduction in the estimated bound is even more significant, passing from 2201~km to 1564~km. Overall, the introduction of regression plays a significant role in the estimate uncertainty reduction, which becomes increasingly important as the sampling time decreases.

Once described the DAIOD algorithm performance, it is now interesting to compare it against other existing methods. Table~\ref{tab:Results_DIM} shows the performance of the \ac{dim} in terms of errors in $a$, $i$, $e$, $\Omega$ and $u$ as a function of $\Delta t_{obs}$ when processing regressed data with $t_s$ = 5~s and $k^{\sigma}_2$. The trend of the errors with the observed arc length recalls the one identified for the DAIOD, but the absolute values are larger. The \ac{dim} provides relatively good results for longer time windows, but its accuracy drops while approaching the short-arcs region. Here the DAIOD method performs significantly better, granting a reduction in the estimate error of at least one order of magnitude. On the other hand, the DIM is generally a factor 2 faster than the DAIOD.

\subsection{Real test cases}
\label{subsec:Real}
This section illustrates the performance of the DAIOD algorithm when processing real data collected by French GRAVES radar. Three different test cases are presented, ordered for increasing values of observed arc length. The first observed object is satellite Cryosat-2. The satellite was observed on 2021-05-02 from 09:08:22 to 09:09:54 UTC, with an overall number of 21~observation instants covering around 1.5$\%$ of the orbital period. The results of the \ac{iod} process are listed in Table~\ref{tab:Results_test1}. The table shows the error $\varepsilon_{DAIOD}$ and bound $b_{DAIOD}$ of the DAIOD method in terms of Keplerian parameters. In addition, the errors of the Doppler integration method $\varepsilon_{DIM}$ are provided. 
\begin{table}[!t]
    \centering
    \small
    \sisetup{round-mode=places,round-precision=4}
    \begin{tabular}{c c c c c c}
    \hline\hline
    & $a$ (km)& $e$ (-)& $i$ (deg)& $\Omega$ (deg)& $u$ (deg)\\ 
    \cline{1-6}
    	$\varepsilon_{DIM}$& 47.822& 0.006& 0.601& 0.454& 0.020\\
		$\varepsilon_{DAIOD}$& \phantom{0}3.312& 0.004& 0.542& 0.423& 0.004\\
		$b_{DAIOD}$& 297.251& 0.404& 1.202& 0.804& 0.335\\
    \hline\hline
    \end{tabular}
    \vspace{0.2cm}
    \caption{DAIOD and \ac{dim} performance comparison for test case 1.}
    \label{tab:Results_test1}
\end{table}
\begin{table}[!t]
    \centering
    \small
    \sisetup{round-mode=places,round-precision=4}
    \begin{tabular}{c c c c c c}
    \hline\hline
    & $a$ (km)& $e$ (-)& $i$ (deg)& $\Omega$ (deg)& $u$ (deg)\\ 
    \cline{1-6}
    	$\varepsilon_{DIM}$& 198.220& 0.026& 0.117& 0.086& 0.011\\
		$\varepsilon_{DAIOD}$& \phantom{0}58.651& 0.004& 0.055& 0.043& 0.059\\
		$b_{DAIOD}$& 116.139& 0.049& 0.319& 0.240& 0.108\\
    \hline\hline
    \end{tabular}
    \vspace{0.2cm}
    \caption{DAIOD and \ac{dim} performance comparison for test case 2.}
    \label{tab:Results_test2}
\end{table}
\begin{table}[!t]
    \centering
    \small
    \sisetup{round-mode=places,round-precision=4}
    \begin{tabular}{c c c c c c}
    \hline\hline
    & $a$ (km)& $e$ (-)& $i$ (deg)& $\Omega$ (deg)& $u$ (deg)\\ 
    \cline{1-6}
    	$\varepsilon_{DIM}$& 111.078& 0.014& 0.251& 0.197& 0.025\\
		$\varepsilon_{DAIOD}$& \phantom{0}34.501& 0.003& 0.066& 0.133& 0.027\\
		$b_{DAIOD}$& 174.054& 0.248& 0.460& 0.355& 0.121\\
    \hline\hline
    \end{tabular}
    \vspace{0.2cm}
    \caption{DAIOD and \ac{dim} performance comparison for test case 3.}
    \label{tab:Results_test3}
\end{table}
The errors are computed with respect to precise orbits whose ephemerides can be retrieved from \texorpdfstring{\protect\footnote{\href{ftp://doris.ign.fr/pub/doris/products/orbits/ssa}{ftp://doris.ign.fr/pub/doris/products/orbits/ssa}}}{}. As can be seen, the DAIOD estimate is systematically more accurate than the \ac{dim} one. Yet, the short duration of the passage unavoidably inflates the estimate bounds, which are two orders of magnitude larger than the actual errors for $a$, $e$, and $u$. This conservative behaviour of the DAIOD method improves when the observed arc length increases. Table~\ref{tab:Results_test2} shows the result of a second test case. The same object was observed by GRAVES on 2021-09-08 from 02:56:18 to 02:58:57 UTC, for an overall number of 61 observation instants covering 2.6 $\%$ of the orbital period. The DAIOD estimate is now less accurate in semi-major axis and more accurate in the angular quantities. On the other hand, all the estimated bounds are smaller than the previous case and enclose well the actual error. The \ac{dim} estimate is still less accurate in all components but $u$. A final test case is shown in Table~\ref{tab:Results_test3}. Radar measurements were collected while observing satellite Sentinel 3A on 2021-09-01 from 21:40:45 to 21:45:01 UTC. Overall, 67 observation instants were collected, covering 4.2 $\%$ of the orbital period. This passage is longer than the previous one, but with a higher sampling time. As a result, no significant reduction in the obtained DAIOD errors and bounds is expected. The results shown in Table~\ref{tab:Results_test3} confirm this guess with both errors and bounds resembling those of Table~\ref{tab:Results_test2}. A similar trend can be noticed for the \ac{dim} too, showing a general lower accuracy than the DAIOD method.

\section{Conclusions}
This paper introduced a novel \ac{iod} algorithm for processing Doppler radar short-arc observations. The method combines \ac{da} and \ac{ads} to obtain an estimate of the state of the space object at the first observation epoch. This estimate is expressed as a set of Taylor polynomials mapping the measurements into the state estimate. As a result, a full description of the orbit set associated with the observations is provided. A detailed analysis of the performance of the method was carried out, and the impact of measurement noise, observed arc length, measurement regression and sampling time was studied. Overall, the DAIOD is more robust than its competitors with both simulated and real observations. More specifically, the analyses show that the method has extremely good convergence properties, regardless of the noise levels and arc lengths. This result is made possible by a proper investigation of the measurements uncertainty set for the first guess computation, which allows the solver to reach success rates above 96~$\%$, even for very short-arc observations with high noise levels. The \ac{iod} solution is then characterised by an accuracy that is strongly dependent on the observation conditions, but is also coupled with a reliable estimate of its bounds. This second result is made possible by the combined use of \ac{da} and \ac{ads}, allowing the DAIOD method to provide a complete Taylor series expansion of the \ac{iod} solution. This aspect is of crucial importance whenever uncertainty-based data association operations shall be performed. Thus, the proposed DAIOD algorithm is a valuable tool whenever radar-based catalogue generation operations are performed.

\section*{Acknowledgment}
This work is founded by CNES through the R-S20/BS-0005-058 grant. The authors thank the French Space Command for providing GRAVES radar measurements to this research. M. Losacco gratefully acknowledges A. Foss\`{a} for the support on the use of CNES' PACE library.

\bibliography{Losacco_DAIOD.bib}

\end{document}